\documentclass[reqno]{amsart}
\usepackage{amsmath}
\usepackage{amssymb}

\usepackage{amsmath, amssymb, amsthm, graphicx, graphics}

\usepackage{color}
\usepackage{ifpdf}
\ifpdf 
  \usepackage[pdftex]{graphicx}
  \DeclareGraphicsExtensions{.pdf,.png,.jpg,.jpeg,.mps}
  \usepackage{pgf}
  \usepackage{tikz}
\else 
  \usepackage{graphicx}
  \DeclareGraphicsExtensions{.eps,.bmp}
  \usepackage{pgf}
  \usepackage{tikz}
  \usepackage{pstricks}
\fi
\usepackage{epic,bez123}
\usepackage{floatflt}
\usepackage{wrapfig}

\newtheorem{lemma}{Lemma}[section]
\newtheorem{theorem}{Theorem}[section]

\theoremstyle{definition}
\newtheorem{definition}{Definition}[section]
\theoremstyle{remark}
\newtheorem{remark}{Remark}[section]

\numberwithin{equation}{section}

\newcommand{\p}{\partial}

\makeatletter
\newcommand{\rmnum}[1]{\mathrm{\romannumeral #1}}
\newcommand{\Rmnum}[1]{\mathrm{\expandafter\@slowromancap\romannumeral#1@}}
\makeatother

\begin{document}
\title[Uniqueness of Transonic Shocks]
{Global Uniqueness of Steady Transonic Shocks in Two--Dimensional
Compressible Euler Flows}

\author{Beixiang Fang}
\author{Li Liu}
\author{Hairong Yuan}

\address{B. Fang:
Department of Mathematics, Shanghai Jiaotong University, Shanghai
200240, China} \email{bxfang@sjtu.edu.cn}

\address{L. Liu:
Department of Mathematics, Shanghai Jiaotong University, Shanghai
200240, China} \email{llbaihe@gmail.com}

\address{H. Yuan (Corresponding Author):
Department of Mathematics, East China Normal University, Shanghai
200241, China} \email{hryuan@math.ecnu.edu.cn;\ \
hairongyuan0110@gmail.com}



\date{\today}

\begin{abstract}
We prove that for the two-dimensional steady complete compressible
Euler system, with given uniform upcoming supersonic flows, the
following three fundamental flow patterns (special solutions) in gas
dynamics involving transonic shocks are all unique in the class of
piecewise $C^1$ smooth functions, under appropriate conditions on
the downstream subsonic flows: $(\rmnum{1})$ the normal transonic
shocks in a straight duct with finite or infinite length, after
fixing a point the shock-front passing through; $(\rmnum{2})$ the
oblique transonic shocks attached to an infinite wedge;
$(\rmnum{3})$ a flat Mach configuration containing one supersonic
shock, two transonic shocks, and a contact discontinuity, after
fixing the point the four discontinuities intersect. These special
solutions are constructed traditionally under the assumption that
they are piecewise constant, and they have played important roles in
the studies of mathematical gas dynamics. Our results show  that the
assumption of piecewise constant can be replaced by some more weaker
assumptions on the downstream subsonic flows, which are sufficient
to uniquely determine these special solutions.

Mathematically, these are uniqueness results on solutions of free
boundary problems of a quasi-linear system of elliptic-hyperbolic
composite-mixed type in bounded or unbounded planar domains, without
any assumptions on smallness. The proof relies on an elliptic system
of pressure $p$ and the tangent of the flow angle $w=v/u$ obtained
by decomposition of the Euler system in Lagrangian coordinates, and
a newly developed method for the $L^{\infty}$ estimate that is
independent of the free boundaries, by combining the maximum
principles of elliptic equations, and careful analysis of shock
polar applied on the (maybe curved) shock-fronts.
\end{abstract}

\maketitle


\section{Introduction}

In discussing steady flows which are supersonic in the entrance
section of a duct and then become subsonic, R. Courant and K.~O.
Friedrichs wrote in their classical monograph {\it Supersonic Flow
and Shock Waves} \cite[p. 372]{CF} that
\begin{quote}
\indent ``$\cdots\cdots$ We know in general that this assumption
(the flow is continuous) is not tenable and that we must consider
flows involving shocks. It is a question of great importance to know
under what circumstances a steady flow involving shocks is uniquely
determined by the boundary conditions at the entrance, and when
further conditions at the exit are appropriate."
\end{quote}
This paper is exactly devoted to showing that, for the
two-dimensional steady complete compressible Euler system, with
given uniform upstream supersonic flows, under certain reasonable
conditions of the downstream of the subsonic flows, in the class of
piecewise $C^1$ functions, the uniqueness of the following three
important flow patterns in gas dynamics:
\begin{itemize}
\item[($\rmnum{1}$)] The transonic normal shock in a two--dimensional straight
duct with finite or infinite length, after fixing a point the
shock-front passing through (see
Figure 1);\\
\item[($\rmnum{2}$)] The transonic oblique shocks
attached to an infinite wedge against the uniform supersonic flow
(see Figure 2);\\
\item[($\rmnum{3}$)] A flat Mach configuration involving three shocks
and a contact discontinuity, where one shock is supersonic, and
another two are transonic, with the contact line separating two
subsonic regions, after fixing a point $O$ where the four
discontinuities intersect (see Figure 3).
\end{itemize}
These flow patterns  have been constructed in \cite[p.311,
pp.332-333]{CF} by using shock polar, under the assumptions that the
shock-fronts (contact line) are straight lines, and the subsonic
flows behind the shock-fronts are also uniform --- i.e., the
solution is piecewise constant (see \S 2 below). The uniqueness
results we obtain in this work indicate that under certain
conditions of the downstream subsonic flows, these assumptions may
be relaxed to that the flow fields are only piecewise continuously
differentiable --- then for given uniform upcoming supersonic flows,
the transonic shock-fronts (and the contact line appeared in the
Mach configuration) must be straight lines and the subsonic flows
behind them must be uniform. From mathematical point of view, since
the steady Euler system for subsonic flows is of elliptic-hyperbolic
composite-mixed type, these are all results on uniqueness of
solutions of free boundary value problems of a nonlinear system of
mixed type, with the transonic shock-fronts and contact line being
the free boundaries.
\begin{figure}[h]
\begin{center}
\includegraphics[scale=0.5]{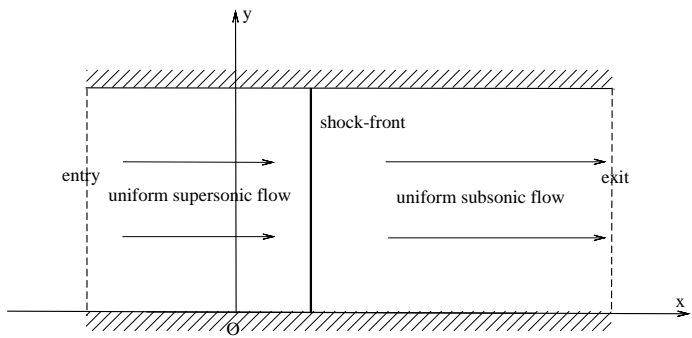}
\end{center}
\caption{{\small Normal transonic shock in a straight duct.}}
\end{figure}
\begin{figure}[h]
\centering
  \setlength{\unitlength}{1bp}%
  \begin{picture}(300, 200)(0,0)
  \put(0,0){\includegraphics[scale=0.5]{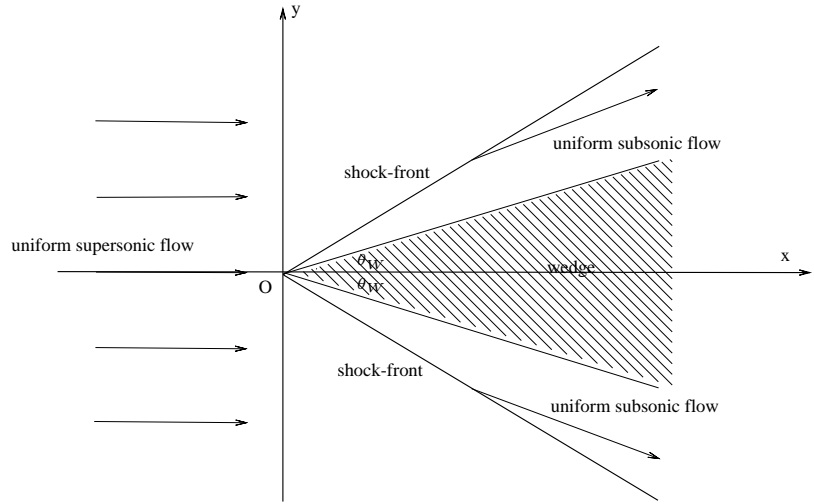}}
  \put(133,92){\fontsize{6}{5}\selectfont ${\tiny\theta_W}$}
  \put(133,83){\fontsize{6}{5}\selectfont ${\tiny\theta_W}$}
  \end{picture}
\caption{{\small  Oblique transonic shocks attached to a slim
infinite symmetric wedge with open angle $2\theta_W$ and zero
attacking angle in uniform supersonic flows.}}
\end{figure}

\begin{figure}[h]
\centering
  \setlength{\unitlength}{1bp}%
  \begin{picture}(300, 200)(0,0)
  \put(0,0){\includegraphics[scale=0.5]{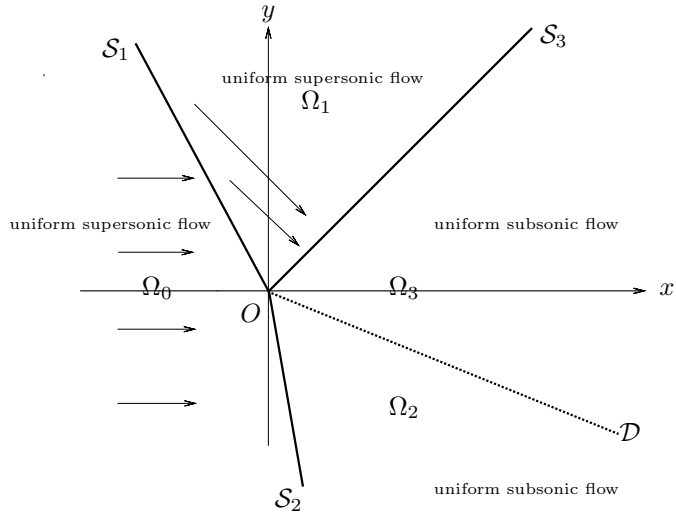}}
  \put(77,65){\fontsize{10}{5}\selectfont ${O}$}
  \put(40,75){\fontsize{10}{5}\selectfont ${\Omega_0}$}
  \put(-10,100){\fontsize{6}{5}\selectfont {uniform supersonic flow}}
  \put(100,145){\fontsize{10}{5}\selectfont ${\Omega_1}$}
  \put(70,155){\fontsize{6}{5}\selectfont {uniform supersonic flow}}
  \put(133,75){\fontsize{10}{5}\selectfont ${\Omega_3}$}
  \put(150,100){\fontsize{6}{5}\selectfont {uniform subsonic flow}}
  \put(133,30){\fontsize{10}{5}\selectfont ${\Omega_2}$}
  \put(150,0){\fontsize{6}{5}\selectfont {uniform subsonic flow}}
   \put(85,180){\fontsize{10}{5}\selectfont $y$}
     \put(235,75){\fontsize{10}{5}\selectfont ${x}$}
     \put(25,165){\fontsize{10}{5}\selectfont ${\mathcal{S}_1}$}
     \put(90,-5){\fontsize{10}{5}\selectfont ${\mathcal{S}_2}$}
     \put(190,170){\fontsize{10}{5}\selectfont ${\mathcal{S}_3}$}
     \put(220,20){\fontsize{10}{5}\selectfont ${\mathcal{D}}$}
  \end{picture}
\caption{\small A flat ``direct" Mach configuration. Here
$\mathcal{S}_1$ is a supersonic shock lying in the quadrant
$\Rmnum{2}$, $\mathcal{S}_2, \mathcal{S}_3$ are transonic shocks
lying in quadrant $\Rmnum{4}, \Rmnum{1}$ respectively, and
$\mathcal{D}$ is a contact line. They are all straight lines. The
flow field is piecewise constant.}
\end{figure}

These special while important flow patterns have been the focus of
research for many years, since the pioneering work on transonic
shocks in ducts by G.-Q. Chen, M. Feldman \cite{cf2003} for the
potential flow equation, and S. Chen \cite{chens1} for the Euler
system. In a series of papers \cite{cf2007,cf2004,cf2003}, G.-Q.
Chen and M. Feldman established the stability of transonic shocks in
finitely long or infinitely long duct with square sections or
arbitrary sections, under various conditions at the exit of the
duct. Z. Xin and H. Yin studied similar problems but employed
different methods \cite{XY2005}. For the Euler system, H. Yuan
showed in \cite{Yu2} various stability and instability results by
using Lagrangian coordinates introduced by S. Chen in \cite{chens2}
and characteristic decomposition technique \cite{chens1}. G.-Q.
Chen, J. Chen and M. Feldman also introduced the method of stream
function in Lagrangian coordinates by which the Euler system can be
written as a single second order equation to study the stability,
local uniqueness and asymptotic behavior of transonic shocks in
infinite ducts for the Euler system. Z. Xin, H. Yin and their
collaborators also studied the problem of stability of transonic
shock in nozzles for the steady Euler system, see
\cite{XYY2009,XY2008}.

For the study of stability of oblique transonic shocks attached to a
wedge in supersonic flow, there are papers of S. Chen, B. Fang
\cite{CHEN-FANG,Fang}, E. H. Kim \cite{Kim}, and H. Yin, C. Zhou
\cite{Yin}. Both the potential flow equation and the complete Euler
system were used as models and the authors showed stability under
suitable downstream conditions.

For the Mach configurations, one has to use the Euler system since
there is a contact discontinuity. S. Chen firstly studied the
stability of a Mach configuration under perturbations of upcoming
supersonic flow, for the steady Euler flow and pesudo-steady Euler
flow \cite{chens2,chen3}. Later on he and B. Fang also showed the
stability of a wave pattern of regular reflection-refraction of
shocks upon an interface \cite{CHEN-FANG}, which is similar to a
Mach configuration.

Comparing to stability, the progress on the study of uniqueness is
rather slow. For potential flows, G.-Q. Chen and H. Yuan \cite{CYu}
proved the transonic shock in a two-dimensional or three-dimensional
straight duct is unique modulo a translation.  L. Liu showed
uniqueness of subsonic potential flow in various domains \cite{L3},
see also \cite{LY2}. By Proposition 3.1 in \cite{ChenJun}, in
two-dimensional case, many of  the uniqueness results on subsonic
flows also hold for the Euler system, since under appropriate
boundary conditions, the flow is in itself irrotational, i.e.,
governed by the potential flow equation. However, to our knowledge,
there is no any result on global uniqueness of transonic shocks for
the complete Euler system before.

In this paper, we are going to establish the global uniqueness, for
the complete Euler system,  of  the three flow patterns we
previously mentioned. It turns out that the method of maximum
principles employed in \cite{CYu} for potential flow also works for
Euler flow, but the point is that one should consider now the
extreme values of pressure $p$ and slope of the flow angle $w=v/u$,
with $u, v$ the velocity component of the flow along $x,y$ axis,
respectively. In \cite{Yu2}, H. Yuan has shown that $p$ and $w$
satisfy a first order elliptic system for subsonic flow. Then as in
\cite{CYu}, one needs to show if there are nontrivial extremes of
$p$ (or $w$) on the shock-front $\mathcal{S}$, there would be
contradictions by a Hopf-type boundary point lemma. It is here we
need to utilize some nice properties of the Rankine-Hugoniot
conditions --- particularly, a rather simple curve called $p-w$
shock polar (see Figure 4), which shows the relation between $p$ and
$w$ behind the shock-front. We recommend \cite{CF,Wjh} for detailed
computations of these relations.

This paper is organized as follows. In \S2, we formulate the above
mentioned physical problems (flow patterns) and state rigorously our
results, Theorem \ref{thm1}--Theorem \ref{thm6}. We also introduce
briefly the $p-w$ shock polar for the comfort of those readers not
familiar with this subject. Then we prove in \S3 the uniqueness of
transonic shocks in straight ducts (Theorem \ref{thm1}--Theorem
\ref{thm2}), and in \S4, the uniqueness of oblique transonic shocks
attached to an infinite wedge (Theorem \ref{thm3}--Theorem
\ref{thm5}). Theorem \ref{thm6}, on the uniqueness of the Mach
configuration studied by S. Chen \cite{chens2}, is established in \S
5.

\section{Shock Polar, Special Solutions and Main Results}

\subsection{The Euler System and $p-w$ Shock Polar}

The complete compressible Euler system governing the steady motion
of polytropic gas flow in two-dimensional space $(x,y)\in
\mathbf{R}^2$ expresses the conservation of mass, momentum and total
energy (cf. \cite[pp.14--23]{CF}):
\begin{eqnarray}
&&\p_x(\rho u)+\p_y(\rho v)=0,\label{e1}\\
&&\p_{x}(\rho u^2+p)+\p_y(\rho uv)=0,\label{e2}\\
&&\p_x(\rho uv)+\p_y(\rho v^2+p)=0,\label{e3}\\
&&\frac 12(u^2+v^2)+\frac{c^2}{\gamma-1}=b_0. \label{e4}
\end{eqnarray}
The last equation is also called as Bernoulli law, with $b_0$ the
Bernoulli constant which is invariant on each flow trajectory, even
across a shock-front. Here $p$, $\rho$, $c$, $u,$ $v$ are the
pressure, mass density, sonic speed, velocity component along $x$
and $y$ axis respectively, of the flow, and $\gamma>1$ is the
adiabatic exponent. For polytropic gas, the state function is
$p=A(s)\rho^\gamma$, with $s$ the entropy, we then have $c=\sqrt{{\p
p}/{\p \rho}}=\gamma p/\rho.$ The Mach number of the flow is defined
by $M=\sqrt{(u^2+v^2)}/c$. When $M<1$ (resp. $M>1$), the flow is
called {\it subsonic flow} (resp. {\it supersonic flow}). It is well
known that the system \eqref{e1}--\eqref{e4} is symmetric hyperbolic
for supersonic flow, while of elliptic-hyperbolic composite-mixed
type for subsonic flow. As in \cite{Yu2}, we introduce $w=v/u$ to
denote the slope of the flow angle \footnote{The flow angle $\theta$
is the angle between the velocity $(u,v)$ and the $x$ axis. Hence
$w=\tan\theta.$} in case $u\ne 0$, which, together with $p$, plays a
key role in our analysis later.

Let $\mathcal{S}=\{(x,y): x=f(y)\}$ be a shock-front --- that is, a
curve where the state of the gas is discontinuous and physical
entropy condition holds (cf. Definition \ref{def101} in \S2.2
below). It is well known that the following Rankine-Hugoniot
conditions should hold across $\mathcal{S}$ \cite[p.11]{Da}:
\begin{eqnarray}
&&[\rho u]-f'(y)[\rho v]=0, \label{rh1}\\
&&[\rho u^2+p]-f'(y)[\rho uv]=0, \label{rh2}\\
&&[\rho uv]-f'(y)[\rho v^2+p]=0. \label{rh3}
\end{eqnarray}
Here $[\cdot]$ denotes the jump of a quantity across $\mathcal{S}$.
For example, if the pressure on the left hand side of $\mathcal{S}$
(upstream) is $p_0$, on the right hand side of $\mathcal{S}$
(downstream) is $p$, then we set $[p]=p-p_0.$ Note that $b_0$ in
\eqref{e4} is constant along a flow trajectory even across the
shock-front, which is the Rankine-Hugoniot condition corresponding
to conservation of energy.
\begin{figure}[h]
\center
\setlength{\unitlength}{1bp}%
  \begin{picture}(300, 200)(40,0)
  \put(0,0){\includegraphics[scale=0.3]{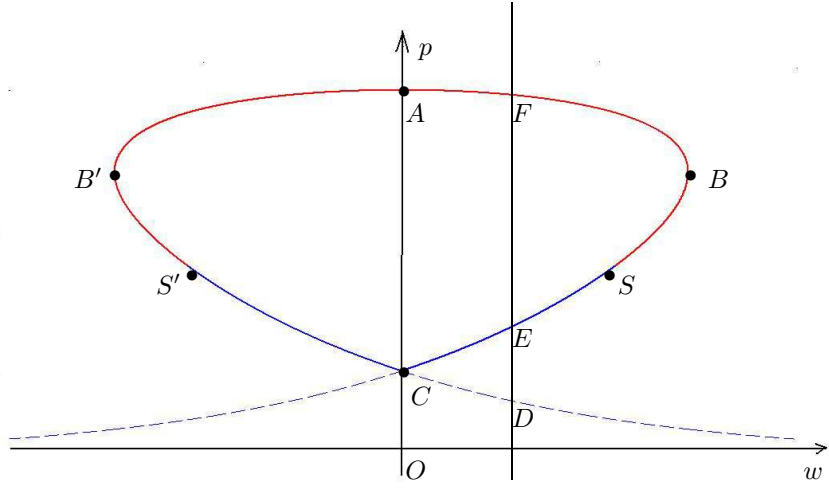}}
  \put(200,20){\fontsize{10}{5}\selectfont ${O}$}
  \put(205,180){\fontsize{10}{5}\selectfont ${p}$}
  \put(200,155){\fontsize{10}{5}\selectfont ${A}$}
   \put(197,164){\fontsize{10}{5}\selectfont ${\bullet}$}
  \put(202,48){\fontsize{10}{5}\selectfont ${C}$}
   \put(197,58){\fontsize{10}{5}\selectfont ${\bullet}$}
  \put(314,130){\fontsize{10}{5}\selectfont ${B}$}
   \put(305,132){\fontsize{10}{5}\selectfont ${\bullet}$}
  \put(75,130){\fontsize{10}{5}\selectfont ${B'}$}
  \put(88,132){\fontsize{10}{5}\selectfont ${\bullet}$}
  \put(280,90){\fontsize{10}{5}\selectfont ${S}$}
  \put(274.5,94.5){\fontsize{10}{5}\selectfont ${\bullet}$}
 \put(106,90){\fontsize{10}{5}\selectfont ${S'}$}
 \put(117,94.5){\fontsize{10}{5}\selectfont ${\bullet}$}
  \put(350,20){\fontsize{10}{5}\selectfont ${w}$}
   \linethickness{0.2pt}
  \put(240,20){\line(0,1){180}}
 \put(240,70){\fontsize{10}{5}\selectfont ${E}$}
 \put(240,40){\fontsize{10}{5}\selectfont ${D}$}
 \put(240,155){\fontsize{10}{5}\selectfont ${F}$}
  \end{picture}

\caption{\small A $p-w$ shock polar. $C=(0,p_0)$ represents the
upcoming supersonic flow; $A=(0,p^+)$ represents the strongest
normal transonic shock; $w$ reaches its maximum at $B=(w_*,p_*)$,
and $S=(w_{sonic}, p_{sonic})$ represents a shock that the state
behind the shock-front is exactly sonic. $B'=(-w_*,p_*)$ and
$S'=(-w_{sonic}, p_{sonic})$ are reflections of $B, S$ with respect
to the $p$ axis. The open arc $\widehat{SAS'}$ represents all
possible transonic shocks. The open arc $\widehat{CS}$ and
$\widehat{CS'}$ represent all possible supersonic shocks. The dashed
line below $C$ corresponds to those states that do not satisfy
physical entropy condition. $E=(w, p_{weak})$ represents the weak
shock; $F=(w, p_{strong})$ represents the strong shock.}
\end{figure}

For a fixed point $Q=(f(y),y)$ on the shock-front, supposing the
uniform state $(p_0>0, u_0>0, w_0\equiv0\equiv v_0, \rho_0>0)$ of
the gas ahead of the shock-front is given \footnote{Therefore $b_0$
in \eqref{e4} is $u_0^2/2+\gamma p_0/((\gamma-1)\rho_0).$}, we
consider the four algebraic equations \eqref{e4}---\eqref{rh3}, with
five independent variables $p(Q), u(Q), w(Q), \rho(Q)$, which are
the state of the flow on the right hand side of $\mathcal{S}$, and
$f'(y).$ In \cite[pp.306--309]{CF}, it has been shown that, the flow
angle $\theta=\arctan w$, horizonal velocity $u(Q)$, and density
$\rho(Q)$, as well as the slope of the shock-front $f'(y)$, can all
be expressed by $p(Q)$. Particularly we have the following important
relation on $p$ and $w$ \cite[p. 347]{CF}:
\begin{eqnarray}\label{pwpolar}
\displaystyle w=\pm \frac{\frac{p}{p_0}-1}{\gamma
M_0^2-(\frac{p}{p_0}-1)}
\sqrt{\frac{\frac{2\gamma}{\gamma+1}(M_0^2-1)-(\frac{p}{p_0}-1)}
{\frac{p}{p_0}+\frac{\gamma-1}{\gamma+1}}}.
\end{eqnarray}
Here $M_0=u_0/c_0$ is the Mach number of the flow ahead of
$\mathcal{S}$, and $c_0=\sqrt{\gamma p_0/\rho_0}.$ Figure 4 shows
the curve represented by this expression in $w-p$ plane, which is
traditionally called a {\it $p-w$ shock polar}. It is easy to see
that this curve is symmetric with respect to the line $w=0$ due to
``$\pm$" in \eqref{pwpolar}.

There are several critical points on the $p-w$ shock polar. $A$ is
the point where $p$ reaches its maximum (denoted by $p^+$); $B$ is
the point where $w$ reaches its maximum $w_*$. The point
$S=(w_{sonic},p_{sonic})$ has the special property that the state of
the gas behind the shock-front is exactly sonic, i.e., $M=1.$ For
all points on the shock polar with $p>p_{sonic}$ (i.e., on the open
arc $\widehat{SAS'}$), the state behind the shock-front they
represent are all subsonic. $B', S'$ are the mirror images of $B,S$
with respect to the line $w=0.$

For convenience of later reference, we also write out the following
relations on the state of the flow behind the shock-front
\cite[p.44, p.48]{Wjh}:
\begin{eqnarray}
u&=&u_0-\frac{p-p_0}{\rho_0u_0},\label{up}\\
\rho&=&\rho_0\frac{(\gamma+1)p+(\gamma-1)p_0}{(\gamma-1)p+(\gamma+1)p_0},\label{rp}\\
\sin\alpha&=&\frac{1}{M_0}\sqrt{\frac{\gamma+1}{2\gamma}\left(\frac{p}{p_0}
+\frac{\gamma-1}{\gamma+1}\right)}, \quad
\cot\alpha=f'(y).\label{fp}
\end{eqnarray}
Here $\alpha$ is the angle between the tangent of the shock-front
and the $x$ axis.

\subsection{The Physical Problems and Main Results}

Now we formulate the physical problems we are concerned with,
construct the special solutions (flow patterns) by using the shock
polar, and then state the results on their uniqueness.

\subsubsection{Problem A: Transonic Shocks in a Finite Duct}
Let $\mathcal{D}=\{(x,y)\in \mathbf{R}^2: -1<x<1, 0<y<1\}$ be a
straight duct with finite length. We set $\Sigma_{-1}=\{(-1,y): 0\le
y\le1\}$ and  $\Sigma_{1}=\{(1,y): 0\le y\le1\}$ as the entry and
exit of the duct, $\Gamma_0=\{(x,0): -1<x<1\}$ and
$\Gamma_1=\{(x,1): -1<x<1\}$ as the lower and upper wall of the
duct, respectively.

On the walls, it is natural to impose the impenetrability or slip
condition $v=0.$ On the entry $\Sigma_{-1}$, we assume the flow is
supersonic and moves in $\mathcal{D}$, so the Euler system
\eqref{e1}--\eqref{e4} is hyperbolic in the positive $x$ direction.
Therefore we impose the initial condition $p=p_0>0, u=u_0>0,
v=v_0=0, \rho=\rho_0>0$, with $p_0, u_0, \rho_0$ all being constants
and satisfying $M_0>1.$ On the exit, one Dirichlet condition on $p$
is given:
\begin{eqnarray}\label{209}
p=p_1,\qquad p_1 \quad\text{is a constant.}
\end{eqnarray}
It is indicated in \cite[p.373, p.385]{CF} that this restriction on
pressure is more physically interesting for nozzle flows. Therefore
we get a boundary value problem of \eqref{e1}---\eqref{e4} in
$\mathcal{D}$. We call this {\it Problem A.} We are interested in
those solutions of  problems like this with the following structure,
i.e., transonic shock solutions:
\begin{definition}\label{def101} Suppose $\mathcal{D}$ is a domain
in $\mathbf{R}^2.$ For a $C^2$ function $x=f(y)$ defined on an
interval $\mathcal{I}$ of $\mathbf{R}$, let
\begin{eqnarray*}
\mathcal{S}&=&\{(f(y), y)\in \mathcal{D} \,:\, y \in \mathcal{I}\},\\
\mathcal{D}^-&=&\{(x,y)\in \mathcal{D}\, :\, x<f(y)\},\\
\mathcal{D}^+&=&\{(x,y)\in \mathcal{D}\,:\, x>f(y)\}.
\end{eqnarray*}
Then $p,u,v,\rho\in L_{loc}^\infty(\mathcal{D})\cap C^1(
{\mathcal{D}^-\cup\mathcal{S}})\cap
C^1({\mathcal{D}^+}\cup\mathcal{S})$ is a {\it transonic shock
solution} of \eqref{e1}--\eqref{e4} if it is supersonic in
$\mathcal{D}^-$ and subsonic in $\mathcal{D}^+$, satisfies equations
\eqref{e1}--\eqref{e4} in $D^-\cup D^+$ and the boundary conditions
pointwise in the classical sense, the Rankine-Hugoniot jump
conditions \eqref{rh1}--\eqref{rh3}, as well as the physical entropy
condition across $\mathcal{S}$: $$[p]>0.$$
\end{definition}

\begin{remark}
By \eqref{rh3} and entropy condition, there holds $f'(y)=\rho
uv/(\rho v^2+[p])$ and the denominator never vanishes. So if
$p,u,v,\rho\in L_{loc}^\infty(\mathcal{D})\cap C^1(
{\mathcal{D}^-}\cup{\mathcal{S}})\cap
C^1({\mathcal{D}^+}\cup\mathcal{S})$, we naturally have $f\in C^2$
in Definition \ref{def101}.
\end{remark}

\begin{remark}
For nozzle flows, we know the transonic shock-front will meet the
walls perpendicularly. So as explained in \cite[p.1348]{Yu2}, at the
corner points where they intersect, the solution will in general be
singular. However, for uniform upcoming supersonic flow and straight
duct, it can be shown that the first order compatibility condition
holds (cf. Lemma 3.3 in \cite[p.19]{XYY2009}), that guarantees the
solution $u,v,p,\rho$ to be $C^1$ at the corners, and hence $C^1$ up
to boundary in $\mathcal{D}$. For oblique shocks attached to a
wedge, in our proof we do not require the solution $u,v,p,\rho$ to
be $C^1$ up to the vertex, where the shock-front meet the wedge
(continuous is enough, see Remark \ref{rem41}). Hence to prove
uniqueness, our assumptions on regularity in Definition \ref{def101}
are reasonable.
\end{remark}

As a simple application of the shock polar, we show existence of
special normal transonic shock solutions to Problem A. For given
uniform upcoming supersonic flow $(p_0>0,u_0>0, w_0=0, \rho_0>0)$, a
$p-w$ shock polar is determined (see Figure 4), and for $w=0,$ we
have a unique $p^+>p_0$ (i.e., the point $A$ in Figure 4)
representing a subsonic state, and hence $u^+,\rho^+$, and
$f'\equiv0$ may be calculated by \eqref{up}\eqref{rp}\eqref{fp} with
$p$ there replaced by $p^+$. Suppose the shock-front passes the
point $(t,0)$ with $-1<t<1$, then we may take the shock-front as
$x\equiv t$, and the state behind it is $U^+=(p^+, u^+, w^+=0,
\rho^+),$ which is uniform and subsonic. It is easy to see that what
we get in this way is a transonic shock solution to Problem A, with
$p_1=p^+$ in \eqref{209}. We denote this solution by $U_b$.

The following theorem is one of the main results of this paper.
\begin{theorem}\label{thm1}
For given supersonic initial data $U_0=(p_0>0,u_0>0,w_0=0,\rho_0>0)$
on $\Sigma_{-1}$, supposing the flow behind the shock-front
satisfies $u>0,$ then Problem A has a transonic shock solution
{\rm(}in the sense of Definition \ref{def101}{\rm)} if and only if
$p_1=p^+$ in \eqref{209}. If this holds, let $U$ be a transonic
shock solution to Problem A with shock-front passing through the
point $(t,0)$, then $U$ must be the special solution $U_b$
constructed above.
\end{theorem}

\subsubsection{Problem $A'$: Transonic Shocks in an Infinite Duct}
This is to consider transonic shocks in an infinitely long straight
duct $\mathcal{D}'=\{(x,y)\in \mathbf{R}^2: -1<x<\infty, 0<y<1\}$.
The boundary conditions are the same as in Problem A, except that
\eqref{209} is dropped. 
We call this boundary value problem of Euler system in the unbounded
domain $\mathcal{D}'$ as {\it Problem $A'$.} Of course $U_b$ is a
solution to this problem. Results like Theorem \ref{thm1} also hold:
\begin{theorem}\label{thm2}
For given supersonic initial data $U_0=(p_0>0, u_0>0,
w_0=0,\rho_0>0)$ on $\Sigma_{-1}$, supposing the flow behind the
shock-front satisfies $u>u_*$ and $M<M^*$ with  fixed numbers
$u_*>0$ and $M^*<1,$ then any transonic shock solution to Problem
$A'$ must be the special solution $U_b$, provided that the
shock-front passes the point $(t,0)\in \overline{\mathcal{D'}}$.
\end{theorem}

Like Theorem 1.3 in \cite[p.568]{CYu}, here we do not need any
asymptotic condition of the flow at infinity, under the rather
strong assumptions that the flow is uniformly subsonic behind the
shock-front, i.e., $M<M_*<1$, and moves uniformly to the  exit at
infinity.

We also note that Theorem \ref{thm1} and Theorem \ref{thm2} imply
particularly the instability of the special normal transonic shocks
under any perturbations of the downstream pressure, as have been
shown in \cite{CYu,Yu2}.

\subsubsection{Problem B: Transonic Shocks Attached to an Infinite Wedge}
Next we consider the transonic oblique shocks attached to an
infinite wedge against a uniform supersonic flow. By symmetry with
respect to the $x$ axis (see Figure 2), we may consider only the
``half-wedge" (a ramp) with an open angle $\theta_W>0$ and zero
attacking angle against the uniform supersonic flow. Let
$\mathcal{D}=\{(x,y)\in\mathbf{R}^2: x>0, y>x\tan\theta_W\}$, and
set $\Sigma_0=\{(x,y): x=0, y\ge0\}$ be the ``entry" where the
uniform supersonic flow  $(p_0>0, u_0>0, w_0=0, \rho_0>0)$
\footnote{Note that this is the initial data to the Euler system on
$\Sigma_0$.} flows from left to right to the ramp
$W=\{(x,y)\in\mathbf{R}^2: x>0, 0<y<x\tan\theta_W\}$ with wall
$\Gamma=\{(x,y): x>0, y=x\tan\theta_W\}$. On the wall we have the
slip condition $w=\tan\theta_W.$ We would like to find transonic
shock in the flow field. We call this {\it Problem B.}

Assuming the transonic shock solution is piecewise constant --- that
is, the shock-front is a straight line issuing from the origin $O$,
which is the vertex of the ramp, and the subsonic flow behind the
shock-front is also uniform, we may obtain special transonic shock
solution to Problem B by the shock polar. Indeed, for
$\tan\theta_W\in(0, w_{sonic})$ (see Figure 4), the line
$w=\tan\theta_W$ meets the $p-w$ shock polar at three different
points $D,E,F$. The point $D$ represents a state that does not
satisfy the physical entropy condition, so we ignore it. The point
$E$ corresponds to a weaker supersonic shock, that is, the flow
behind the shock-front is also supersonic. This supersonic shock has
been intensively studied (see, for example,
\cite{chenli,czz,LIYU,Zhangyongqian2003} and references therein),
and is not of our concern. The point $F$ represents a transonic
shock which we are interested in, and we may obtain the slope of the
shock-front, the density and velocity behind the shock-front by
\eqref{up}--\eqref{fp}. We call this solution the {\it strong
transonic shock $U_s$.} It always corresponds to the point
$(w=\tan\theta_W\in(0,w_*], p=p_{strong})$ on the upper part of the
shock polar (i.e., on the closed arc $\widehat{AB}$). For
$w=\tan\theta_W\in (w_{sonic}, w_*),$ both the weak shock and strong
shock are transonic shocks, and for $w=w_*$, they coincide. We will
denote the weak transonic shock by $U_w$; it corresponds to the
point $(w=\tan\theta_W\in(w_{sonic},w_*), p=p_{weak}<p_{strong})$ on
the lower part of the shock polar (i.e., on the open arc
$\widehat{BS}$). For $w> w_*$, detached bow shock will appear, and
for $w=w_{sonic}$, the weak shock is sonic. We will not consider
these two cumbersome cases in this paper.

For the special strong transonic shock, we have the following
uniqueness result.
\begin{theorem}\label{thm3}
For given supersonic initial data $U_0=(p_0>0,u_0>0,w_0=0,\rho_0>0)$
on $\Sigma_0$, supposing the flow behind the shock-front satisfies
$u>0$ and the following asymptotic condition:
\begin{eqnarray}\label{212}
w\rightarrow\tan\theta_W  \qquad \text{as}\ \
(x,y)\rightarrow\infty,
\end{eqnarray}
here $\tan\theta_W\in(0,w_*],$ then any transonic shock solution to
Problem B with the restriction at the vertex of the ramp
\begin{eqnarray}\label{213}
w=\tan\theta_W, \quad p=p_{strong}\qquad \text{at}\ \ O
\end{eqnarray}
must be the strong transonic shock solution $U_s.$
\end{theorem}

Under a further rather strong asymptotic condition on the
shock-front at infinity, namely,
\begin{eqnarray}\label{214}
p<p_{\infty}\qquad \text{as}\ \  (x,y)\rightarrow\infty\ \
\text{on}\ \ S
\end{eqnarray}
with $p_{\infty}$ a constant less than $p_*,$ we also have
uniqueness of the weak transonic shock solution.
\begin{theorem}\label{thm4}
For given supersonic initial data $U_0=(p_0>0,u_0>0,w_0=0,\rho_0>0)$
on $\Sigma_0$, supposing the flow behind the shock-front satisfies
$u>0,$ and the asymptotic conditions \eqref{212}, \eqref{214}, then
any transonic shock solution to Problem B with the restriction at
the vertex of the ramp
\begin{eqnarray}\label{215}
w=\tan\theta_W, \quad p=p_{weak}\qquad   \text{at}\ \ O
\end{eqnarray}
must be the weak transonic shock solution $U_w.$ Here
$\tan\theta_W\in (w_{sonic},w_*).$
\end{theorem}

Here we note that by shock polar, one can always determine that at
$O$ either \eqref{213} or \eqref{215} holds. Therefore our
requirements in the Theorems are natural.

We may replace \eqref{212} in Theorem \ref{thm3}  by
\begin{eqnarray}\label{217}
p\rightarrow p_{strong}  \qquad \text{as}\ \ (x,y)\rightarrow\infty
\end{eqnarray}
to guarantee  uniqueness of $U_s$, and replace \eqref{212} in
Theorem \ref{thm4} by
\begin{eqnarray}\label{218}
p\rightarrow p_{weak}  \qquad \text{as}\ \ (x,y)\rightarrow\infty,
\end{eqnarray}
as well as dropping \eqref{214} to show uniqueness of $U_w$.
\begin{theorem}\label{thm5}
For given supersonic initial data $U_0=(p_0>0,u_0>0,w_0=0,\rho_0>0)$
on $\Sigma_0$, supposing the flow behind the shock-front satisfies
$u>0,$ and the asymptotic condition \eqref{217} {\rm(}resp.
\eqref{218}{\rm)}, then any transonic shock solution to Problem B
with the restriction \eqref{213} {\rm(}resp. \eqref{215}{\rm)} at
the vertex of the ramp must be the strong {\rm(}resp. weak{\rm)}
transonic shock solution $U_s$ {\rm(}resp. $U_w${\rm)}.
\end{theorem}

We recall that, in discussing the shocks attached to a wedge, R.
Courant and K. O. Friedrichs had written
\cite[pp.317-318]{CF}:``Quite aside from the question of stability,
the problem of determining which of the possible shocks occurs
cannot be formulated and answered without taking the boundary
conditions at infinity into account. $\cdots\cdots$ If the pressure
prescribed there is below an appropriate limit, the weak shock
occurs in the corner. If, however, the pressure at the downstream
end is sufficiently high, a strong shock may be needed for
adjustment. Under appropriate circumstances this strong shock may
begin just in the corner and thus, of the two possibilities
mentioned, the one giving a strong shock may actually occur.
$\cdots\cdots$ All statements made here are conjectures so far.
While there is little doubt that they are in general correct, they
should be supported, if possible, by detailed theoretical
investigation." Out results above, may be considered as a part of
such investigation.

\subsubsection{Problem C: Uniqueness of a Mach Configuration}
Mach configuration is one of the fundamental flow patterns in gas
dynamics which is closely connected to shock reflection-diffraction
phenomena. As explained by R. Courant and K. O. Friedrichs in
\cite[pp.332-333]{CF}: ``Three shocks separating three zones of
different continuous states are impossible.  $\cdots\cdots$
Configurations involving three shock fronts must therefore involve
an additional discontinuity. The simplest assumption, in agreement
with many observations, is that there occurs in addition to the
three shock fronts a single contact discontinuity line." Such a wave
structure is called a {\it Mach configuration} (cf.
\cite[p.2]{chens2}). Figure 3 shows a possible Mach configuration,
where $\mathcal{S}_1$ is an incident shock separating the two
supersonic regions $\Omega_0$ and $\Omega_1$; $\mathcal{S}_3$ is
called the reflected shock from the contact discontinuity line
$\mathcal{D}$; $\mathcal{S}_2$ is called a Mach front. The flow in
$\Omega_2$ and $\Omega_3$ are both subsonic. So $\mathcal{S}_1$ is a
supersonic shock and $\mathcal{S}_2, \mathcal{S}_3$ are transonic
shocks. Following terminologies in \cite{chens2}, we call
$\Omega_0\cup\mathcal{S}_1\cup\Omega_1$ the supersonic part of the
Mach configuration, and
$\mathcal{S}_3\cup\Omega_3\cup\mathcal{D}\cup\Omega_2\cup
\mathcal{S}_2$ its subsonic part.

For the case that all the discontinuities are straight lines and the
flow field is piece-wise constant, the existence of such a ``flat"
Mach configuration showed in Figure 3 can be proved by using shock
polar. In Figure 5, the point $C$ represents the upstream supersonic
flow in $\Omega_0$. $I_1$ corresponds to the flow state in
$\Omega_1$. By the shock polar $R$ roots at $I_1$, we can determine
all the possible states in $\Omega_3$; and by  the  polar roots at
$C$, we may determine all the possible states in $\Omega_2$. Since
$w$ and $p$ should be continuous across the contact line
$\mathcal{D}$, the point
\begin{equation}\label{pm}
I_{2,3}=(w_m, p_m)
\end{equation}
where the two loops intersect represents the $w$ and $p$ in
$\Omega_2, \Omega_3$.  Hence by \eqref{up}--\eqref{fp}, the slope of
the lines $\mathcal{S}_1$, $\mathcal{S}_2$, $\mathcal{S}_3$,
$\mathcal{D}$ and $u,v$, as well as  $\rho$, in $\Omega_1, \Omega_2,
\Omega_3,$ can all be uniquely determined.  We recommend
\cite[pp.346-350]{CF} or \cite[pp.22-24]{Ben} for a detailed
discussion.
\begin{figure}[h]
\centering
  \setlength{\unitlength}{1bp}%
  \begin{picture}(300, 200)(40,0)
  \put(0,0){\includegraphics[scale=0.3]{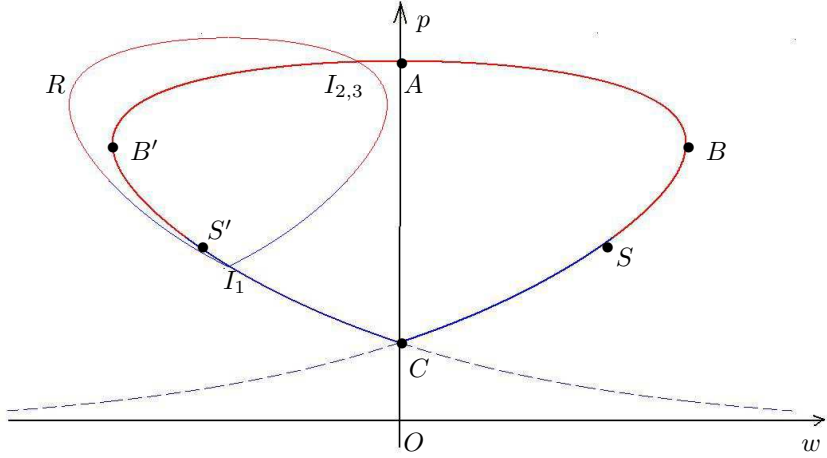}}
  \put(200,20){\fontsize{10}{5}\selectfont ${O}$}
  \put(205,180){\fontsize{10}{5}\selectfont ${p}$}
  \put(200,155){\fontsize{10}{5}\selectfont ${A}$}
   \put(197,164){\fontsize{10}{5}\selectfont ${\bullet}$}
  \put(202,48){\fontsize{10}{5}\selectfont ${C}$}
   \put(197,58){\fontsize{10}{5}\selectfont ${\bullet}$}
  \put(314,130){\fontsize{10}{5}\selectfont ${B}$}
   \put(305,132){\fontsize{10}{5}\selectfont ${\bullet}$}
  \put(97,130){\fontsize{10}{5}\selectfont ${B'}$}
  \put(88,132){\fontsize{10}{5}\selectfont ${\bullet}$}
  \put(280,90){\fontsize{10}{5}\selectfont ${S}$}
  \put(274.5,94.5){\fontsize{10}{5}\selectfont ${\bullet}$}
 \put(126,100){\fontsize{10}{5}\selectfont ${S'}$}
 \put(122,94.5){\fontsize{10}{5}\selectfont ${\bullet}$}
  \put(350,20){\fontsize{10}{5}\selectfont ${w}$}
  \put(132,81){\fontsize{10}{5}\selectfont ${I_1}$}
  \put(170,155){\fontsize{10}{5}\selectfont ${I_{2,3}}$}
  \put(65,155){\fontsize{10}{5}\selectfont ${R}$}
  \end{picture}

\caption{\small Existence of a flat Mach configuration by using the
$p-w$ shock polar.}
\end{figure}

It arises naturally the question of uniqueness: under what
conditions, can we show the flat Mach configuration constructed
above is the only possible one in the class of piecewise $C^1$
functions?  We have the  follows.
\begin{theorem}\label{thm6}
Suppose the supersonic part, i.e., the uniform upcoming supersonic
flow in $\Omega_0$ and a straight supersonic shock-front
$\mathcal{S}_1$, are given (see Figure 6), and $u>0$, as well as the
following asymptotic condition in the subsonic region $\Omega_2,
\Omega_3$:
\begin{eqnarray}\label{pa}
p\rightarrow p_m \quad \text{as}& (x,y)\rightarrow\infty,
\end{eqnarray}
where $p_m$ is determined in \eqref{pm}. For definiteness, we also
fix the point $O\in\mathcal{S}_1$ where it meets the possibly curved
discontinuities $\mathcal{S}_2, \mathcal{S}_3, \mathcal{D}.$ Then
any weak entropy solution $(p,w, u,\rho) \in
L_{loc}^\infty(\mathbf{R}^2)\cap C^1(\overline{\Omega_0})\cap
C^1(\overline{\Omega_1}) \cap
C_{loc}^1(\overline{\Omega_2}\backslash\{O\})\cap
C_{loc}^1(\overline{\Omega_3}\backslash\{O\})$ of the Euler system
must be the flat Mach configuration. Here
$C_{loc}^1(\overline{\Omega}\backslash\{O\})$ is the set of those
functions continuous in
 $\overline{\Omega}$ and continuously differentiable in any compact
subset of $\overline{\Omega}\backslash\{O\}$. \eqref{pa} may be
replaced by
\begin{eqnarray}\label{wa}
w\rightarrow p_m \quad \text{as}& \Omega_2\cup\Omega_3\ni
(x,y)\rightarrow\infty,
\end{eqnarray}
and the same conclusion holds.
\end{theorem}
\begin{figure}[h]
\centering
  \setlength{\unitlength}{1bp}%
  \begin{picture}(300, 200)(0,0)
  \put(0,0){\includegraphics[scale=0.5]{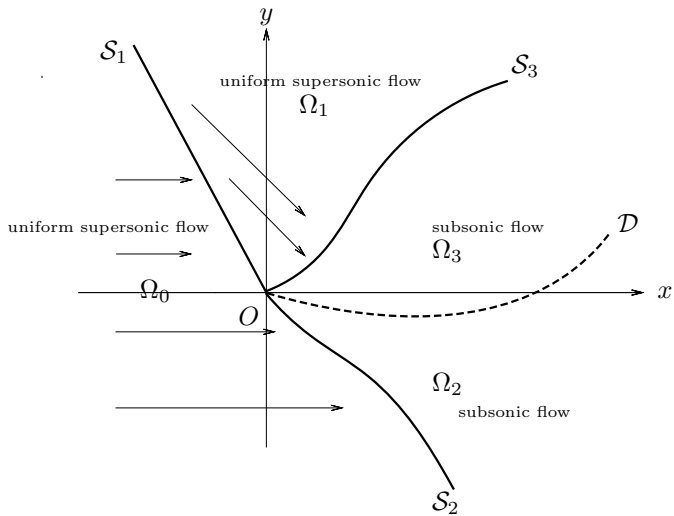}}
  \put(77,65){\fontsize{10}{5}\selectfont ${O}$}
  \put(40,75){\fontsize{10}{5}\selectfont ${\Omega_0}$}
  \put(100,145){\fontsize{10}{5}\selectfont ${\Omega_1}$}
  \put(150,90){\fontsize{10}{5}\selectfont ${\Omega_3}$}
  \put(150,40){\fontsize{10}{5}\selectfont ${\Omega_2}$}
   \put(85,180){\fontsize{10}{5}\selectfont $y$}
     \put(235,75){\fontsize{10}{5}\selectfont $x$}
     \put(25,165){\fontsize{10}{5}\selectfont ${\mathcal{S}_1}$}
     \put(150,-5){\fontsize{10}{5}\selectfont ${\mathcal{S}_2}$}
     \put(180,160){\fontsize{10}{5}\selectfont ${\mathcal{S}_3}$}
     \put(220,100){\fontsize{10}{5}\selectfont ${\mathcal{D}}$}
     \put(-10,100){\fontsize{6}{5}\selectfont {uniform supersonic flow}}
      \put(70,155){\fontsize{6}{5}\selectfont {uniform supersonic flow}}
      \put(160,30){\fontsize{6}{5}\selectfont {subsonic flow}}
      \put(150,100){\fontsize{6}{5}\selectfont {subsonic flow}}
  \end{picture}
\caption{\small The possible curved Mach configuration in Euler
coordinates $(x,y)$.}
\end{figure}

\smallskip
For the four Problems stated above, by the theory of classical
solutions of hyperbolic systems \cite{LIYU}, it is easy to show that
the supersonic flow ahead of the transonic shock-fronts must be
uniform (cf. \cite[p.216]{chenchenf} for a sketch). Particularly, we
may prescribe supersonic initial data on $\{x=0\}$ in Problem C (see
Figure 6). To prove these uniqueness results, one then only needs to
show the state behind the transonic shock must be uniform, and hence
by shock polar, the shock-front must be straight. We will show how
these hold and prove Theorem \ref{thm1}, Theorem \ref{thm2} in \S 3,
Theorem \ref{thm3}-- Theorem \ref{thm5} in \S 4, and Theorem
\ref{thm6} in \S 5, by using a rather unified and newly developed
approach that combining various maximum principles and analysis of
the shock polar.

\section{Uniqueness of Transonic Shocks in Ducts}

\subsection{Reduction of the Euler System}

Supposing
\begin{eqnarray}\label{add3}
\rho u>0 \quad\text{and}\quad 0<M<1,
\end{eqnarray}
the Euler system \eqref{e1}--\eqref{e4} may be written in the
Lagrangian coordinates $(\xi,\eta)$ as (see
\cite[pp.1351--1356]{Yu2})
\begin{eqnarray}
&\p_\xi p+\lambda_R\p_\eta p-\beta_1\p_\eta w=0,\label{a01}\\
&\p_\xi w+\beta_2\p_\eta p+\lambda_R\p_\eta w=0,\label{a02}
\end{eqnarray}
together with the invariance of entropy along flow trajectories
\begin{equation}
\p_\xi\left(\frac{p}{\rho^\gamma}\right)=0,\label{a03}
\end{equation}
and the Bernoulli law
\begin{equation}\label{a04}
\frac12u^2(1+w^2)+\frac{c^2}{\gamma-1}=b_0.
\end{equation}
Here $\xi=x$, and $\eta$ measures the mass flux along two flow
trajectories. The flow trajectory in Lagrangian coordinates is the
line $\eta=\mathrm{constant}$. The coefficients in \eqref{a01} and
\eqref{a02} are
\begin{eqnarray}\label{coe}
\begin{cases}\displaystyle
\lambda_R=\frac{\rho c^2 uw}{u^2-c^2},\\
\displaystyle \beta_1=\frac{-\rho^2c^2u^3}{u^2-c^2}>0,\\
\displaystyle \beta_2=\frac{(M^2-1)c^2}{u(u^2-c^2)}>0.
\end{cases}
\end{eqnarray}
These may be obtained from (4.8)(4.19) in \cite[pp.1355--1356]{Yu2},
by taking $\psi=\psi'=0$ there. Note our assumption $0<M<1$ in
\eqref{add3} implies $u\ne c$, which is required in (4.21) in
\cite[p.1356]{Yu2}. By (3.11) or (3.12) in \cite[p.1352]{Yu2} and
our assumptions in Definition \ref{def101}, the transformation
$(x,y)\mapsto(\xi,\eta)$ is $C^2$ and invertible, provided that
$\rho u>0$. It is easy to check that for subsonic flows, \eqref{a01}
and \eqref{a02} constitute an elliptic system, which is the elliptic
part of the Euler system. Equations \eqref{a03}\eqref{a04} are the
hyperbolic part of the Euler system.

The equations \eqref{a01} \eqref{a02} may also be expressed as
\begin{eqnarray}\label{a05}
Dw=WDp,
\end{eqnarray}
with $Dw=(\p_\xi w, \p_\eta w)^T$ the gradient of a function $w$,
and $W$ the matrix
\begin{eqnarray}\label{a06}
W=\frac{-1}{\beta_1}\left(
    \begin{array}{cc}
      \lambda_R & \beta_1\beta_2+\lambda_R^2 \\
      -1 & -\lambda_R \\
    \end{array}
  \right).
\end{eqnarray}
Note that $\det W=\beta_2/\beta_1=(1-M^2)/(\rho^2u^4)>0,$ so $W$ is
invertible. By acting $(\p_\eta, -\p_\xi)$ from left to \eqref{a05},
we get a second order elliptic equation of $p$ in divergence form:
\begin{eqnarray}\label{a07}
\p_{\xi}\left(\frac{\p_\xi
p}{\beta_1}+\frac{\lambda_R}{\beta_1}\p_\eta
p\right)+\p_\eta\left(\frac{\lambda_R}{\beta_1}\p_\xi p
+(\beta_2+\frac{\lambda_R^2}{\beta_1})\p_\eta p\right)=0.
\end{eqnarray}
By our assumptions in Definition \ref{def101}, the coefficients
\footnote{They are $1/\beta_1,\ \lambda_R/\beta_1$,
$\beta_2+\lambda_R^2/\beta_1$.} are all $C^1$ functions
\footnote{Hence bounded in any fixed compact domain.}, and the
equation is uniformly elliptic in any fixed compact domain.
Therefore the strong maximum principle for weak solutions of
elliptic equations in divergence form  is available, see Theorem
8.19 in \cite[p. 189]{GT}. We also need the following boundary point
lemma due to R. Finn and D. Gilbarg (see Lemma 7 in
\cite[p.31]{FG}). \footnote{Although the formulation here is a
little different from that in \cite{FG}, this lemma is obviously a
direct consequence of that.}
\begin{lemma}\label{Hopf}
Let the coefficients of the system
\begin{eqnarray}\label{lu}
-v_y=au_x+b'u_y,\qquad v_x=b''u_x+cu_y,\\
\triangle=4ac-(b'+b'')^2\ge\delta>0\label{3111}\quad\quad
\end{eqnarray}
be H\"older continuous in the closure of a simply connected domain
$\Omega\subset \mathbf{R}^2$ with $C^{1,\alpha}$ boundary
{\rm(}$0<\alpha<1${\rm)}, and $u, v\in C^1(\overline{\Omega})$ is a
solution to \eqref{lu}. Suppose there is a $Q\in\p\Omega$ satisfying
$u(x)>u(Q)$ {\rm($u(x)<u(Q)$)} for all $x\in \Omega,$ then there
holds
\begin{eqnarray}\label{H1}
\frac{\p u}{\p n}(Q)>0, \qquad{\rm(}\frac{\p u}{\p n}(Q)<0,{\rm)}
\end{eqnarray}
where $n$ is the inward drawn normal of the curve $\p\Omega$ at $Q$.
\end{lemma}

It is clear that the inner normal $n$ in \eqref{H1} may be replaced
by any vector $l$ at $Q$ with $l\bullet n>0,$ since along the
tangent $\tau$ of $\p \Omega$ we have $\frac{\p u}{\p\tau}(Q)=0$. In
this paper we always use ``$\bullet$" to denote the scalar product
of vectors.

\subsection{Finitely Long Duct: Proof of Theorem \ref{thm1}}

Under the assumptions in Definition \ref{def101} and Theorem
\ref{thm1}, we show \eqref{add3} holds. We directly have $u>0$ and
$M<1$. By \eqref{pwpolar}, \eqref{rp} and \eqref{a03}, the entropy
$A(s)=p/\rho^\gamma$ has the estimate
\begin{eqnarray}\label{add4}
0<s_*<A(s)<s^* \qquad \text{in}\ \ \mathcal{D}^+
\end{eqnarray}
with two numbers $s_*, s^*$ depending only on the upcoming
supersonic flow $p_0, u_0, \rho_0$ and the adiabatic  exponent
$\gamma>1.$ Hence from $u>0$ and $M<1$ we get $0<u^2<c^2=\gamma
A(s)\rho^{\gamma-1}$, therefore $\rho>0$ in $\mathcal{D}^+.$ One
then checks that for \eqref{a05},
$\triangle=4\beta_2/\beta_1=4(1-M^2)/(\rho^2u^4)$ is bounded away
from zero in a compact subset of $\overline{\mathcal{D}^+}$. So we
may employ the results in \S 3.1 for the following analysis.

The duct, in the Lagrangian coordinates, is now the rectangle
$\{(\xi,\eta): -1<\xi<1, 0<\eta<\eta_0\},$ with $\eta_0=\rho_0u_0>0$
(see Figure 7). At the exit $\Sigma_1=\{(1,\eta):
0\le\eta\le\eta_0\}$ we still have the Dirichlet condition for
pressure:
\begin{equation}\label{a08}
p=p_1,
\end{equation}
with $p_1$ a constant (cf. \eqref{209}). On the upper and lower
walls $\Gamma_1=\{(\xi,\eta_0): -1<\xi<1\}, \Gamma_0=\{(\xi,0):
-1<\xi<1\},$ there should hold $w=0,$ hence $\p_\xi w=0,$ and by
\eqref{a05}, this is a Neumann condition of $p$:
\begin{eqnarray}\label{a09}
\p_\eta p=0,
\end{eqnarray}
since now  we have  $\lambda_R=0$ and $\beta_1\beta_2>0.$
\begin{figure}[h]
\centering
  \setlength{\unitlength}{1.00bp}%
  \begin{picture}(250, 130)(0,0)
  \put(0,0){\includegraphics[scale=0.5]{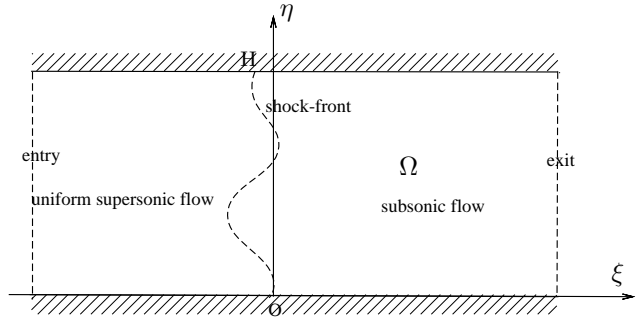}}
  \put(105,120){\fontsize{10}{17.07}\selectfont $\eta$}
\put(90,101){\fontsize{8}{10}\selectfont H}
  \put(230,20){\fontsize{10}{17.07}\selectfont $\xi$}
  \put(150,60){\fontsize{10}{17.07}\selectfont $\Omega$}
  \end{picture}\caption{{\small The duct and transonic shock in $(\xi,\eta)$
coordinates.}}
\end{figure}

Without loss of generality, we suppose the shock-front
$\mathcal{S}$, with the equation $\xi=\psi(\eta)\in C^2[0,\eta_0]$
in Lagrangian coordinates,  passing through the origin $O$ (i.e.,
$\psi(0)=0$ or $f(0)=0$). We do not need write out the
Rankine-Hugoniot conditions in Lagrangian coordinates here because
the relations \eqref{pwpolar}---\eqref{rp} do not depend on
coordinate system (these Rankine-Hugoniot conditions are available
in \cite[p.1358]{Yu2}). Since $w=0$ on the walls, the shock-front
always meets the walls perpendicularly (i.e.,
$\psi'(0)=0=\psi'(\eta_0))$, even in the $(\xi,\eta)$ coordinates
\footnote{See (6.5) in \cite[p.1358]{Yu2}, that is,
$\psi'(\eta)=[uw]/[p].$}.

Let $\Omega$ be the domain bounded by the shock-front, walls and
exit, where the flow is subsonic.  Now if $p$ is a constant in
$\Omega$, that is, $p\equiv p_1$ according to \eqref{a08}, then by
\eqref{a05}, $w$ is a constant, and by boundary conditions on walls,
we have $w\equiv 0$. So from the shock polar, we can infer that
$p_1=p^+.$ If $p_1\ne p^+$, then we get a contradiction and either
the solution is not constant or  $C^1(\bar{\Omega})$ solution does
not exist. We further may solve on $\mathcal{S}$ that $p=p^+, u=u^+,
\rho=\rho^+$ and $f'=0$ by \eqref{pwpolar}---\eqref{fp}, and by the
hyperbolic part \eqref{a03}--\eqref{a04}, obtain that $\rho=\rho^+,
u=u^+$ in $\Omega$. Hence if $p$ is a constant in $\Omega,$ then
$p_1=p^+$ and the solution is exactly the special solution $U_b$.
The uniqueness is therefore proved.

Now we show that if $p$ is not a constant in $\Omega$, then there
will occur contradictions by using maximum principles. In the
following, we always assume that $p$ is not a constant in $\Omega$.

{\sc Step 1.}\ First, by \eqref{a07} and the strong maximum
principle (Theorem 8.19 in \cite[p.189]{GT}), $p_{max}$, $p_{min}$,
the maximum and minimum of $p$ in $\overline{\Omega}$, respectively,
can only be attained on the boundary
$\p\Omega=\Sigma_1\cup(\Gamma_0\cap\overline{\Omega})\cup
(\Gamma_1\cap\overline{\Omega})\cup{\mathcal{S}}.$ Since $p_1$ is a
constant, $p_{max}$ and $p_{min}$ can not be attained simultaneously
on $\Sigma_1.$ By  Lemma \ref{Hopf} and \eqref{a09}, none of them
can be attained on the walls. \footnote{Here the walls are the line
segments  $\Gamma_i\cap\overline{\Omega}$ ($i=0,1$) excluding their
extreme points, hence they are open. Note that, however, by our
definition, $\Sigma_1$ and $\mathcal{S}$ are closed.} So at least
one of them should be attained on ${\mathcal{S}}$.

{\sc Step 2.}\ Now we suppose that $p_{max}$ is attained on
${\mathcal{S}}$. By shock polar, we can see the maximum is $p^+$
(corresponding to $w=0$, that is the point $A$ in Figure 4) and is
certainly attained at $O=(0,0)$ and $H=(\psi(\eta_0),\eta_0)$. Let
$\tau=(\psi'(\eta), 1)$ and $n=(1,-\psi'(\eta))$ be respectively the
tangential and inner normal vector of $\mathcal{S}$ at
$(\psi(\eta),\eta)$, where $p$ attains its maximum. By \eqref{a05},
we have
\begin{eqnarray}\label{a15}
\frac{\p w}{\p\tau}=\tau Dw=(\tau W)Dp.
\end{eqnarray}
But
\begin{eqnarray}\label{a16}
(\tau W)\bullet n=\tau W
n^T=\frac{1}{\beta_1}((\lambda_R\psi'-1)^2+\beta_1\beta_2\psi'^2)>0,
\end{eqnarray}
by Lemma \ref{Hopf} we get \begin{equation}\label{add1}(\tau
W)Dp<0.\end{equation}

Note that at the corner points $O, H,$ we cannot employ  Lemma
\ref{Hopf} directly. However, since both the corner angles are
$\pi/2$, we may use a reflection argument to overcome this
difficulty. Taking $O$ as an example. We first reflect $\mathcal{S},
\Omega$ with respect to $\eta=0$ and obtain its image $\mathcal{S}',
\Omega'.$ Since the corner angle is $\pi/2,$ the curve
$\mathcal{S}\cup \mathcal{S}'$ is still $C^2.$ Let $B_r(O)$ be an
open disk centered at $O$ with a very small radius $r$, we set
$\mathcal{O}=B_r(O)\cap(\Omega\cup\Omega'\cup\Gamma_0)$. Then we use
even reflection to extend $p$ (both the solution $p$ and those
appeared in the coefficients through \eqref{a05}) to $\mathcal{O}$,
and by Neumann condition \eqref{a09}, $p\in C^1(\mathcal{O}).$
Similarly we use even reflection to extend $\rho, u$, odd reflection
to extend $v,w$ to $\mathcal{O}.$  Hence in $\mathcal{O}$ equation
\eqref{a07} still holds and the coefficients are at least Lipschtiz
continuous. We then use Lemma \ref{Hopf} at $O\in \mathcal{O}$ to
obtain \eqref{add1}.

Hence whenever $p$ reaches its maximum, we have $\frac{\p w}{\p
\tau}<0$ by \eqref{a15} and \eqref{add1}. Thus from $w(O)=0,$ we may
infer that $w$ is always negative on ${\mathcal{S}}\backslash\{O\}$.
This is impossible since $w(H)=0.$

{\sc Step 3.}\ Next we consider the case that only the minimum of
$p$ is attained on $\mathcal{S}$. If $p_{min}=p^+,$ then by shock
polar, we see that on $\mathcal{S}$ there should always hold $p=p^+$
and $w=0,$ hence $\psi'=0.$ We then get a contradiction from
\eqref{a15} and \eqref{a16}: the left hand side of \eqref{a15} is
zero, while the right hand side should be positive due to  Lemma
\ref{Hopf}.

{\sc Step 4.}\ We then consider the case that  $p_{min}<p^+.$  First
we show $p_{min}\ne p_*$. Indeed, if $p_{min}=p_*$, by shock polar
there are only two possibilities: the points corresponding to $B$ or
$B'$ in Figure 4, where $p=p_*$, and $w$ attains its maximum
$w_*$(at $B$), or minimum $-w_*$ (at $B'$). By \eqref{a16}, we still
see the left hand side of \eqref{a15} is zero, while the right hand
side is positive, a contradiction!

{\sc Step 5.}\ Then we conclude that $p_{min}\in (p_{sonic},
p_*)\cup (p_*, p^+).$ We remark that $p_{min}$ should be larger than
$p_{sonic}$ by Definition \ref{def101}, since for those points on
shock polar below $S$ and $S'$, the flow is supersonic, and sonic at
$S,S'$.

Considering $p$ as a function of $w$ for $p\in(p_{sonic}, p_*)$ and
$p\in(p_*, p^+)$ respectively, we see that $p'(w)\ne0.$  Now suppose
$p$ attains its minimum $p_{min}$ at a point $P$ on $\mathcal{S}$,
then at $P$ we have
\begin{eqnarray}\label{add7}
0=\tau Dp=\tau (p'(w)Dw)=p'(w)\tau Dw=p'(w)\tau W Dp
\end{eqnarray}
due to \eqref{a05}. Therefore $\tau W Dp=0,$ a contradiction to
Lemma \ref{Hopf} by \eqref{a16}, which yields $\tau W Dp>0$. This
finishes the proof of Theorem \ref{thm1}.

\begin{remark}
We may  replace \eqref{209} or \eqref{a08},  the Dirichlet condition
of $p$, by $w=0$ on the exit $\Sigma_1.$ The proof is similar to
those presented below in \S 4.1, and we omit the details.
\end{remark}

\subsection{Infinitely Long Duct: Proof of Theorem \ref{thm2}}

For infinitely long duct, we set $\Omega'$ as the subsonic region
behind the shock-front $\mathcal{S}$:
\begin{eqnarray}
\Omega'=\{(\xi,\eta): \xi>\psi(\eta), 0<\eta<\eta_0\}.
\end{eqnarray}
We first show that under the assumptions in Theorem \ref{thm2},
there holds
\begin{eqnarray}\label{add5}
||p,\rho, u,v||_{C(\overline{\Omega'})}\le C_0,
\end{eqnarray}
with $C_0$ depending only on $p_0,u_0,\rho_0$ and $\gamma$. Indeed,
by Bernoulli law, the bounds of $u, v$ and $c$ is
$\max\{\sqrt{2},\sqrt{\gamma-1}\}\sqrt{b_0}.$ But
$\rho=\left({c^2}/{(\gamma A(s))}\right)^{1/(\gamma-1)}$, by
\eqref{add4}, $\rho$ is bounded. Hence due to $p=c^2\rho/\gamma$,
$p$ is also bounded.

We claim the coefficients in equation \eqref{a07} is bounded, with a
bounds determined by $p_0, \rho_0, u_0$ and $\gamma, u_*, M^*$. To
prove this, from \eqref{add4}, we have a lower bound on $\rho$:
\[
\rho=\left(\frac{c^2}{\gamma
A(s)}\right)^{\frac{1}{\gamma-1}}>\rho_*:=\left(\frac{u_*^2}{\gamma
s^*}\right)^{\frac{1}{\gamma-1}}.
\]
Then by \eqref{coe}, it is easy to check that
\begin{eqnarray*}
&&|\lambda_R|\le\frac{\rho |v|}{1-M^2}\le \frac{C_0^2}{1-{M^*}^2},\\
&&\rho_*^2u_*^3<\beta_1\leq\frac{C_0^5}{1-{M^*}^2},\\
&&\frac {1}{C_0}\leq\beta_2\le \frac{1}{u_*(1-{M^*}^2)},
\end{eqnarray*}
and the claim is proved.

Finally, since $M<M^*<1,$ it is straightforward to check that the
equation \eqref{a07} is uniformly elliptic in $\Omega'$, and
\eqref{3111} holds in $\Omega'$.

We now state an extreme principle of  bounded solutions to general
linear elliptic equations in the unbounded domain $\Omega'$ with
Neumann conditions on the lateral boundaries.
\begin{lemma}\label{lem2}
Let $Lu=\p_{\xi}(a^{11}\p_\xi u+a^{12}\p_\eta
u)+\p_{\eta}(a^{21}\p_\xi u+a^{22}\p_\eta u)=0$  be a linear
uniformly elliptic equation with bounded coefficients $a^{ij}\in
C^1_{loc}(\overline{\Omega'}).$
\footnote{$C^1_{loc}(\overline{\Omega'})$ is the set of those
functions which are $C^1$ in any compact subset of
$\overline{\Omega'}$.} Suppose $u\in C(\overline{\Omega'})\cap
C_{loc}^1(\overline{\Omega'})$ is a bounded weak solution
\footnote{A $W_{loc}^{1,2}(\Omega')$ function $u$ is a weak solution
to a mixed boundary value problem if it satisfies the Dirichlet
condition on $\mathcal{S}$ in the sense of trace, the Neumann
condition on $\Gamma_0, \Gamma_1$ as stated in Lemma \ref{lem3}
below.} to this equation that satisfies the homogeneous Neumann
conditions on the lateral boundaries $\Gamma_0'=\{(\xi,\eta):
\xi>\psi(0), \eta=0\}$ and $\Gamma_1'=\{(\xi,\eta):
\xi>\psi(\eta_0), \eta=\eta_0\}.$ Then either
$u_{inf}=\inf_{\overline{\Omega'}}u$ or
$u_{sup}=\sup_{\overline{\Omega'}}u$ is attained on $\mathcal{S}.$
\end{lemma}
By this lemma,  we see  at least one of the extremes, $p_{inf}$ and
$p_{sup}$, which are respectively the infimum and supremum of $p$ in
$\overline{\Omega'}$, is attained on ${\mathcal{S}}$. The analysis
in \S3.2 then shows the uniqueness of the transonic shock solution,
hence Theorem \ref{thm2} is proved. \footnote{It is worth to note
that we need not assume $u,v,p,\rho\in C^1(\overline{\Omega'})$ in
the proof of Theorem \ref{thm2}. The assumption $u,v,p,\rho\in
C_{loc}^1{(\overline{\Omega'})}$ on regularity is enough.}

To prove  Lemma \ref{lem2}, we need an up to boundary Harnack
inequality for elliptic equations with Neumann boundary conditions:
\begin{lemma}\label{lem3}
Let $\Omega\subset\mathbf{R}^n$ be a domain with part of its
boundary $T\subset\p\Omega$ lying in a plane, and
$\sum_{i,j=1}^n\p_i(a^{ij}\p_ju)=0$ a linear equation with
coefficients $a^{ij}$ satisfying
\begin{eqnarray}
\sum_{i,j=1}^na^{ij}\xi_i\xi_j\ge \lambda|\xi|^2& &\text{for any}\ \
(\xi_1,\cdots,\xi_n)\in\mathbf{R}^n,\label{321}\\
\sum_{i,j=1}^n|a^{ij}|\le \Lambda & &\text{in}\ \ \Omega\label{322}
\end{eqnarray}
for two fixed positive constants $\lambda, \Lambda.$ Suppose $u\in
W_{loc}^{1,2}(\Omega)$ is a continuous nonnegative weak solution
satisfying the homogeneous Neumann condition on $T$, that is, for
any $\varphi\in C_c^\infty(\mathbf{R}^n)$ whose support
$\mathrm{supp}\,\varphi$ is disjoint with $\p\Omega\backslash T$,
there holds
\begin{eqnarray}
\int_{\mathrm{supp}\,\varphi\cap (\Omega\cup T)}\sum_{i,j=1}^n
a^{ij}\p_ju\p_i\varphi\,dx=0.
\end{eqnarray}
Then for any compact subset $K$ of $\Omega\cup T$, there holds
\begin{eqnarray}\label{harnack}
\max_K u\le C\min_K u,
\end{eqnarray}
with the positive constant $C$ relying only on $K, \Omega,n,
\Lambda/\lambda.$
\end{lemma}
\begin{proof}
Let $B_r(Q)$ be the open ball centered at $Q\in \mathbf{R}^n$ with
radius $r.$ For $K=\overline{B_r(Q)}\subset\subset \Omega,$
\eqref{harnack} is the classical interior Harnack inequality
(Theorem 8.20 in \cite[p.199]{GT}).

We then consider the case that $Q\in T.$ Without loss of generality,
we may suppose $T$ lying in the hyperplane
$\mathbf{R}^n_0:=\{x=(x_1,\cdots,x_n)\in\mathbf{R}^n: x_n=0\}.$ By
even reflections of $a^{ij}$($1\le i,j\le n-1$ or $i=j=n$) and $u$
with respect to $\mathbf{R}^n_0$, and odd reflections of
$a^{ij}$($i=n$ or $j=n$ but $i\ne j$) with respect to
$\mathbf{R}^n_0$, they are all extended from
$B_{2r}^+(Q):=B_{2r}(Q)\cap\{x=(x_1,\cdots,x_n)\in\mathbf{R}^n:
x_n>0\}\subset \Omega$ to the whole ball $B_{2r}(Q)$ for a small
$r$. We also get an elliptic equation satisfying \eqref{321}
\eqref{322} in $B_{2r}(Q)$. It is straightforward to check that $u$
is still a nonnegative weak solution (in the sense of (8.2) in
\cite[p.177]{GT}) to this equation. Hence by the interior Harnack
inequality we also have $\max_{\overline{B_r(Q)}} u\le
C\min_{\overline{B_r(Q)}} u,$ and \eqref{harnack} follows with
$K=\overline{B_r^+(Q)}.$

For $K$  any compact subset of $\Omega\cup T$, we may cover it by
finite balls (or half balls with center on $T$), and \eqref{harnack}
holds by an argument as that in the proof of Theorem 2.5 in
\cite[p.16]{GT}.
\end{proof}

\noindent{\it Proof of Lemma \ref{lem2}.} If $u$ is a constant, then
obviously Lemma \ref{lem2} holds. In the following we always assume
that $u$ is not a constant. We prove by contradiction, that is, we
assume that
\begin{eqnarray}\label{assumption}
u_{inf}<\min_\mathcal{S} u,\qquad
u_{sup}>\max_\mathcal{S}u.\end{eqnarray} The argument below is
similar to that of proof of Theorem 1.2 in \cite[p.570]{CYu}.

Set $\Omega'_L=\Omega'\cap\{\xi<L\}$ and
$\Sigma_L=\Omega'\cap\{\xi=L\}.$ By \eqref{assumption}, and applying
Lemma \ref{Hopf}, as well as the strong maximum principle in
$\Omega'_L$, we infer that there is an increasing  sequence
$\{L_k\}$ tends to infinity, and there are points $P_k, Q_k\in
\overline{\Sigma_{L_k}},$ such that
\begin{eqnarray}
u(P_k)&=&\max_{\overline{\Omega'_{L_k}}} u\nearrow u_{sup},\\
u(Q_k)&=&\min_{\overline{\Omega'_{L_k}}} u\searrow  u_{inf}
\end{eqnarray}
as $k\rightarrow\infty.$

We may choose $L_1$  rather large. Then define $U_k=\{(\xi,\eta)\in
\Omega': L_k-2<\xi<L_k+2, 0\le \eta\le\eta_0\},$ and
$V_k=\{(\xi,\eta)\in \Omega': L_k-1\le\xi\le L_k+1, 0\le
\eta\le\eta_0\}.$ By  translation along $\xi$: $$\xi=\xi'+L_k,\quad
\eta=\eta',$$ we may translate $U_k$ to $U_0=\{(\xi',\eta'):
-2<\xi'<2, 0\le \eta'\le\eta_0\}$, and $V_k$ to $V_0=\{(\xi',\eta'):
-1\le\xi'\le1, 0\le \eta'\le\eta_0\},$ and the equation $Lu=0$ in
$U_k$ are transformed to $L'_ku_k=0$ in $U_0$, with
$u_k(\xi',\eta')=u(\xi'+L_k, \eta').$

Now let $v_k=u_{sup}-u_k,$ which is positive  and also  a weak
solution of $L'_kv=0$ in $U_0$. Applying Lemma \ref{lem3} to $v_k$,
by the strong maximum principle, we have
\begin{eqnarray}
u_{sup}-u(Q_k)<\max_{V_0} v_k\le C\min_{V_0}v_k<C(u_{sup}-u(P_k)).
\end{eqnarray}
It is crucial to note that by our assumptions in Lemma \ref{lem2},
the constant $C$ here actually does not depend on $k$. Then taking
$k\rightarrow\infty,$ we have
\begin{eqnarray}
0<u_{sup}-u_{inf}<C\times 0=0,
\end{eqnarray}
a contradiction desired. This finishes the proof of Lemma
\ref{lem2}.

\section{Uniqueness of  Transonic Shocks Attached
to a Wedge }

\subsection{Uniqueness of Strong Transonic Shock: Proof of Theorem \ref{thm3}}

Now we consider the oblique transonic shocks attached to an infinite
wedge against uniform supersonic flow (see Figure 2).
\begin{figure}[h]
\centering
   \setlength{\unitlength}{1bp}%
  \begin{picture}(400, 180)(0,0)
  \put(0,0){\includegraphics[scale=0.5]{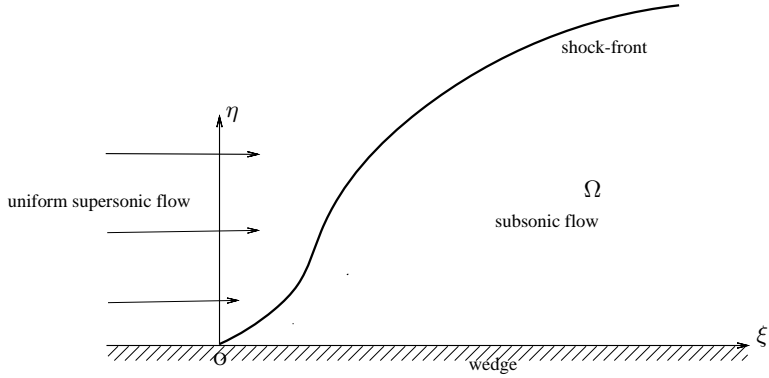}}
  \put(285,15){\fontsize{9}{9}\selectfont $\xi$}
  \put(85,100){\fontsize{9}{9}\selectfont $\eta$}
  \put(220,70){\fontsize{9}{9}\selectfont $\Omega$}
  \end{picture}
  \caption{
\small The oblique transonic shock attached to a wedge in Lagrangian
 coordinates.}
\end{figure}

In Lagrangian coordinates, the wall of the wedge is
$\Gamma=\{(\xi,\eta): \xi>0, \eta=0\},$ see Figure 8. On $\Gamma$
the slip condition is still
\begin{equation}\label{add2}w=\tan\theta_W.\end{equation}
As in \S 3.2, let the equation of the shock-front $\mathcal{S}$ be
$\xi=\psi(\eta)$ with $\eta\ge0$ and $\psi(0)=0,$ and
$\Omega=\{(\xi,\eta): \eta>0, \psi(\eta)<\xi<\infty\}$. Condition
\eqref{213} is now
\begin{eqnarray}\label{a18}
w=\tan\theta_W,\quad p=p_{strong}
\end{eqnarray} which holds at the
origin $O$, and \eqref{212} is replaced by
\begin{eqnarray}\label{a19}
w\rightarrow\tan\theta_W \qquad \text{as}\ \ (\xi,
\eta)\rightarrow\infty.
\end{eqnarray}
By the latter requirement and shock polar, we infer that on
$\mathcal{S}$, there holds
\begin{eqnarray}\label{a20}
p\rightarrow p_{strong} \ \ \text{or}\ \ p\rightarrow p_{weak}
\qquad \text{as}\quad (\xi, \eta)\rightarrow\infty.
\end{eqnarray}

We derive a second order elliptic equation of $w$. Indeed, by
\eqref{a05}, we have
\begin{eqnarray}\label{a21}
Dp=W^{-1}Dw, \qquad W^{-1}=-\frac{\beta_1}{\beta_2}W,
\end{eqnarray}
then $w$ satisfies
\begin{eqnarray}\label{a22}
\p_{\xi}\left(\frac{\p_\xi
w}{\beta_2}+\frac{\lambda_R}{\beta_2}\p_\eta
w\right)+\p_\eta\left(\frac{\lambda_R}{\beta_2}\p_\xi w
+(\beta_1+\frac{\lambda_R^2}{\beta_2})\p_\eta w\right)=0,
\end{eqnarray}
and \eqref{a16} is now
\begin{eqnarray}\label{a23}
\tau W^{-1}n^T=-\frac{\beta_1}{\beta_2}\tau Wn^T<0.
\end{eqnarray}
One may check that Lemma \ref{Hopf} holds for \eqref{a21}, and the
strong maximum principle holds for $\eqref{a22}.$

If $w$ is a constant, i.e., $w\equiv \tan\theta_W$ in $\Omega$, then
by \eqref{a18} and shock polar, we infer that there should hold
$p=p_{strong}$ on $\mathcal{S}$ and $\mathcal{S}$ is the special
strong transonic shock-front by \eqref{fp}. From \eqref{up} and
\eqref{rp}, $u$ and $\rho$ can also be calculated on $\mathcal{S}$.
Then using the hyperbolic part \eqref{a03}\eqref{a04}, we show the
solution must be $U_s$, and Theorem \ref{thm3} is proved.

In the following we show there are contradictions if $w$ is not a
constant. Suppose now that $w$ is not identical to $\tan\theta_W$ in
$\Omega$. Then by Lemma \ref{Hopf} and the strong maximum principle,
either the supremum of $w$ in $\overline{\Omega}$, $w_{sup}$, or the
infimum of $w$ in $\overline{\Omega},$ $w_{inf}$, must be attained
on $\mathcal{S}$. (By \eqref{add2} and \eqref{a19}, $w_{sup}$ is
indeed the maximum of $w$, while $w_{inf}$ is the minimum of $w$.)

We consider firstly the case that $\tan\theta_W\in (0, w_*).$

{\sc Step 1.}\ Suppose $w_{sup}$ is attained at a point $P$ on
$\mathcal{S}$. If $\tan\theta_W<w_{sup}< w_*$, then $P\ne O.$ Note
that by shock polar, we may regard $w$ as a function of $p$,
$w=w(p)$, which shares $w'(p)\ne0$ for $w\ne \pm w_*.$ Now
considering the directional derivative of $w$ along $\mathcal{S}$ at
$P$, $\frac{\p w}{\p\tau}(P)$, which vanishes, we have
\begin{eqnarray}
0=\frac{\p w}{\p\tau}(P)=\tau Dw(P)=w'(p)\tau Dp(P)=w'(p)\tau
W^{-1}Dw(P).
\end{eqnarray}
Hence $\tau W^{-1}Dw(P)=0.$ However, by \eqref{a23} and Lemma
\ref{Hopf}, one gets $\tau W^{-1}Dw(P)>0$, a contradiction!

{\sc Step 2.}\ Secondly consider the case  $w_{sup}=w_*,$ whence
$p(P)=p_*, P\ne O$. Then by \eqref{a21} and  Lemma \ref{Hopf}, we
have
\begin{eqnarray}\label{a25}
\tau Dp=\tau W^{-1}Dw>0.
\end{eqnarray}
This means $p$ is strictly increasing once it attains $p_*$.
However, let $P$ be the first point on $\mathcal{S}$ where $w$
attains its maximum as $\eta$ varies from $0$ to $\infty$. Note that
$p_{strong}>p_*$, so $p$ decreases near the left side of $P$, we
should have $\tau Dp(P)\le 0,$ a contradiction to the assumption of
$p\in C^1(\mathcal{S})$.

{\sc Step 3.}\ We then conclude that  $w_{sup}=\tan\theta_W,$ which
is attained at $O.$ So $w_{inf}<\tan\theta_W,$ and definitely
$w_{inf}$ is attained on $\mathcal{S}\backslash\{O\}$. If
$w_{inf}\ne -w_*$ and $w_{inf}\ne 0$, then as in {\sc Step 1} above
we can show a contradiction since $w'(p)\ne 0$.

{\sc Step 4.}\ Now for $w_{inf}=-w_*,$ as that of \eqref{a25}, we
may get $\tau Dp=\tau W^{-1}Dw<0.$ Therefore, whenever $w$ attains
its minimum $-w_*$, $p$ is strictly decreasing (hence less than
$p_*$) and $w$ should run into $(-w_*, -w_{sonic})$ (the arc
$\widehat{B'S'}$ of shock polar), this contradicts to our asymptotic
condition \eqref{a19}: $w\rightarrow\tan\theta_W>0.$

{\sc Step 5.}\ For the case $w_{inf}=0,$ where $p$ attained its
maximum on $\mathcal{S}\backslash\{O\}$, we get $0=\tau Dp=\tau
W^{-1}Dw<0$, a contradiction.

Secondly, for the case that $\tan\theta_W=w_*$, we have
$w(O)=w_{sup}=w_*$ by \eqref{a18} and \eqref{a19}. This fact,
together with \eqref{add2}, \eqref{a19} and strong maximum
principle, implies  $w_{inf}<w_*$ should be attained on
$\mathcal{S}\backslash\{O\}.$ We then repeat {\sc step 3 -- step 5}
above to reach contradictions.

This finishes proof of Theorem \ref{thm3}.

\subsection{Uniqueness of Weak Transonic Shock: Proof of Theorem \ref{thm4}}

Next we consider the uniqueness of weak transonic shock. Now at the
vertex of the ramp $O$, we have
\begin{eqnarray}
w=\tan\theta_W,\quad p=p_{weak}
\end{eqnarray}
by \eqref{215}. The condition at infinity is \eqref{a19} plus
\begin{eqnarray}\label{a27}
p< p_{\infty}\qquad \text{as}\ \  \eta\rightarrow\infty\ \
\text{on}\ \ \mathcal{S}
\end{eqnarray}
with $p_{\infty}$ a constant less than $p_*.$

As in \S 4.1, to prove uniqueness, we only need to  show that there
are contradictions if $w$ is not a constant, and then either its
supremum $w_{sup}$ or infimum $w_{inf}$ is attained on
$\mathcal{S}$.

{\sc Step 1.}\ Now consider the case $w_{sup}$ is attained on
$\mathcal{S}$. If $\tan\theta_W<w_{sup}<w_*$, the proof is the same
as in {\sc Step 1} in \S 4.1. For $w_{sup}=w_*$, by \eqref{a25},
once $w$ reaches its maximum, $p$ should increase (larger than
$p_*$). This would be a contradiction to our assumption \eqref{a27}
and \eqref{a20}.

{\sc Step 2.} For the case that $w_{inf}<\tan\theta_W$ is attained
on $\mathcal{S}$, if $w_{inf}$ lies on the arc $\widehat{BAS'}$, by
continuity, $w$ should first attain its maximum $w_*$ (note now $w$
runs counterclockwise on the shock polar from $(\tan\theta_W,
p_{weak})$ to $B'$ on the arc $\widehat{BAS'}$). There is a
contradiction as we showed above in {\sc Step 1}; if $w_{inf}$ lies
on the arc $\widehat{SB}$, just repeat {\sc Step 1} in \S 4.1 for a
contradiction.

Theorem \ref{thm4} is also proved.

\subsection{Proof of Theorem \ref{thm5}}

Finally, we show Theorem \ref{thm5} holds. The arguments are very
similar to that in \S 3.2, so we just outline the main ideas.

By \eqref{a05} and \eqref{add2}, one gets an oblique derivative
condition of $p$ on $W$:
\begin{eqnarray}\label{412}
\lambda_R\p_\xi p+(\beta_1\beta_2+\lambda_R^2)\p_\eta p=0.
\end{eqnarray}
The asymptotic conditions \eqref{217} \eqref{218} are respectively
\begin{eqnarray}
&p\rightarrow p_{strong}  \qquad &\text{as}\ \
(\xi,\eta)\rightarrow\infty, \label{413}\\
&p\rightarrow p_{weak} \qquad &\text{as}\ \
(\xi,\eta)\rightarrow\infty.\label{414}
\end{eqnarray}
To show uniqueness, as before, it is sufficient to show
contradictions if $p$ is not constant. In the following we always
assume $p$ is not identical to $p_{strong}$ (resp. $p_{weak}$). Then
by \eqref{412}--\eqref{414} and strong maximum principles applied to
equation \eqref{a07}, either $p_{sup}=\sup_{\overline{\Omega}}p$ or
$p_{sup}=\inf_{\overline{\Omega}}p$ is attained on $\mathcal{S}.$

\subsubsection{Uniqueness of Strong Transonic Shock}

 {\sc Step 1.} Consider the case that $p_{sup}>p_{strong}$ is attained at
 a point $(\psi(\eta_*),\eta_*)\in\mathcal{S}$. If
 $p_{sup}=p^+$, then as of Step 2 in \S 3.2, $w$ should be
 negative on $\mathcal{S}$ for $\eta>\eta_*$. Note that \eqref{413}
 implies particularly $w\rightarrow\pm\tan\theta_W.$ So we get
 $w\rightarrow-\tan\theta_W<0.$ Hence $\psi'(\eta)=uw/[p]$
 \footnote{See footnote 5. Note that $[p]>0$ and $u>0$ are both bounded
 away from zero on the shock-front which can be checked by using shock polar.} has a negative
 upper bound for large $\eta$, which means $\mathcal{S}$ will
 eventually intersect with the boundary $\{(\xi,\eta): \xi=0, \eta>0\}$
 where hyperbolic initial data is prescribed, that is not allowed.

{\sc Step 2.} If $p^+>p_{sup}>p_{strong}$, then $p_{sup}$ lies on
the arc $\widehat{AB}$ of $p-w$ shock polar, and by \eqref{add7}
there is a contradiction. So $p_{sup}=p_{strong}.$ Since $p$ is not
constant, we infer $p_{inf}<p_{strong}$ and is attained on
$\mathcal{S}\backslash\{O\}$. Just repeat {\sc Step 4--Step 5} in
\S3.2 for contradictions. This finishes the proof of Theorem
\ref{thm5} for the case of strong transonic shocks. Note that the
arguments here also work for the case $\tan\theta_W=w_*.$

\subsubsection{Uniqueness of Weak Transonic Shock}
 Consider as well the case that $p_{sup}>p_{weak}$ is attained at a
point $(\psi(\eta_*),\eta_*)\in\mathcal{S}$. If
 $p_{sup}=p^+$, then the argument is the same as in {\sc Step 1} in \S 4.3.1.
If $p_{sup}=p_*$, refer {\sc Step 4} in \S 3.2 for a contradiction.
If $p_{sup}\in (p_{weak}, p_*)\cup(p_*, p^+),$ refer {\sc Step 5}
there. Hence we conclude that $p_{sup}=p_{weak},$ and
$p_{inf}<p_{weak}$ is attained on $\mathcal{S}\backslash\{O\}.$
Therefore $p_{inf}\in (p_{sonic}, p_*)$ and refer also {\sc Step 5}
in \S 3.2 for a contradiction. This finishes proof of Theorem
\ref{thm5} for the case of weak transonic shock.

\begin{remark}\label{rem41}
In the above  proofs we do not need to take care of the corner
points $O$, so we may only require that $p,\rho,w,u\in
C_{loc}^1(\overline{\Omega}\backslash\{O\})\cap
C_{loc}(\overline{\Omega}).$\footnote{$C_{loc}(\overline{\Omega})$
is the set of those functions that are continuous in any compact
subset of $\overline{\Omega}$.}
\end{remark}

\section{Uniqueness of a Mach Configuration}

Since the contact line $\mathcal{D}$ is also a trajectory of the
flow, by virtue of Lagrangian coordinates, it may be transformed to
the $\xi$ axis, see Figure 9. In the following we use superscript
$\cdot^{(i)}$, such as  $u^{(i)},p^{(i)},w^{(i)},\rho^{(i)}$ to
denote quantities in $\Omega_i$, with $i=2,3$. $\mathcal{D}$ is a
contact line means
\begin{equation}\label{501}
w^{(2)}=w^{(3)}, \ \  p^{(2)}=p^{(3)} \qquad\text{on}\ \  \eta=0.
\end{equation}
That is, $w$ and $p$ are continuous in $\Omega_2\cup\Omega_3$,
although others, such as, $\rho, u, v, s$, are not. The supersonic
initial data is prescribed on the $\eta$ axis.
\begin{figure}[h]
\centering
  \setlength{\unitlength}{1bp}%
  \begin{picture}(300, 200)(0,0)
  \put(0,0){\includegraphics[scale=0.5]{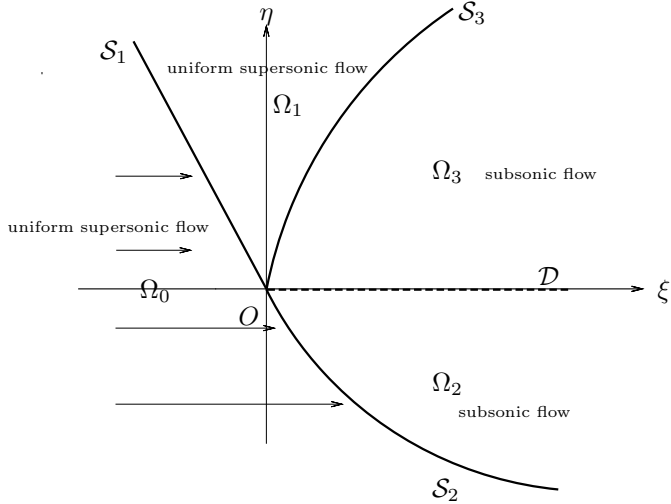}}
  \put(77,65){\fontsize{10}{5}\selectfont ${O}$}
  \put(40,75){\fontsize{10}{5}\selectfont ${\Omega_0}$}
  \put(90,145){\fontsize{10}{5}\selectfont ${\Omega_1}$}
  \put(150,120){\fontsize{10}{5}\selectfont ${\Omega_3}$}
  \put(150,40){\fontsize{10}{5}\selectfont ${\Omega_2}$}
   \put(85,180){\fontsize{10}{5}\selectfont $\eta$}
     \put(235,75){\fontsize{10}{5}\selectfont $\xi$}
     \put(25,165){\fontsize{10}{5}\selectfont ${\mathcal{S}_1}$}
     \put(150,0){\fontsize{10}{5}\selectfont ${\mathcal{S}_2}$}
     \put(160,180){\fontsize{10}{5}\selectfont ${\mathcal{S}_3}$}
     \put(190,80){\fontsize{10}{5}\selectfont ${\mathcal{D}}$}
     \put(-10,100){\fontsize{6}{5}\selectfont {uniform supersonic flow}}
      \put(50,160){\fontsize{6}{5}\selectfont {uniform supersonic flow}}
      \put(160,30){\fontsize{6}{5}\selectfont {subsonic flow}}
      \put(170,120){\fontsize{6}{5}\selectfont {subsonic flow}}
  \end{picture}
\caption{\small The possible curved Mach configuration in Lagrangian
coordinates $(\xi,\eta)$.}
\end{figure}

We note that by the analysis of shock polar in \S 2.2.4, we can
uniquely solve $p^{(2)}(O), w^{(2)}(O), u^{(2)}(O), \rho^{(2)}(O)$,
$p^{(3)}(O), w^{(3)}(O), u^{(3)}(O), \rho^{(3)}(O),$ as well as the
slope of the shock-fronts $\mathcal{S}_2, \mathcal{S}_3$ at $O$.
Particularly, we have (cf. \eqref{pm})
\begin{eqnarray}\label{502}
p^{(2)}(O)=p^{(3)}(O)=p_m,\quad w^{(2)}(O)=w^{(3)}(O)=w_m.
\end{eqnarray}

In the subsonic region $\Omega_i, i=2,3,$ the Euler system reads
\begin{eqnarray}
&\p_\xi p^{(i)}+\lambda_R^{(i)}\p_\eta p^{(i)}-\beta_1^{(i)}\p_\eta w^{(i)}=0,\label{503}\\
&\p_\xi w^{(i)}+\beta_2^{(i)}\p_\eta p^{(i)}+\lambda_R^{(i)}\p_\eta w^{(i)}=0,\label{504}\\
&\p_\xi\left(\frac{p^{(i)}}{(\rho^{(i)})^\gamma}\right)=0,\label{505}\\
&\label{506}
\frac12(u^{(i)})^2(1+(w^{(i)})^2)+\frac{(c^{(i)})^2}{\gamma-1}=b_0.
\end{eqnarray}
Here the coefficients are
\begin{eqnarray}\label{507}
\begin{cases}\displaystyle
\lambda_R^{(i)}=\frac{\rho^{(i)} (c^{(i)})^2
u^{(i)}w^{(i)}}{(u^{(i)})^2-{(c^{(i)})}^2},\\
\displaystyle
\beta_1^{(i)}=\frac{-(\rho^{(i)})^2(c^{(i)})^2(u^{(i)})^3}
{(u^{(i)})^2-(c^{(i)})^2}>0,\\
\displaystyle
\beta_2^{(i)}=\frac{((M^{(i)})^2-1)(c^{(i)})^2}{u^{(i)}
((u^{(i)})^2-(c^{(i)})^2)}>0.
\end{cases}
\end{eqnarray}

Now we begin to prove Theorem \ref{thm6}. It is divided into several
steps. We consider only the case of the condition \eqref{pa}. The
case of \eqref{wa} may be proved similarly by adopting arguments in
\S 4.1 and then combining the following reasoning, so we omit it.

{\sc Step 1.}  If both $p^{(2)}$ and $p^{(3)}$ are constants in
$\Omega_2$, $\Omega_3$ respectively, then we see the subsonic flow
in $\Omega_2, \Omega_3$ are uniform, and by shock polar, the
transonic shock-fronts are straight lines, and then according to
analysis in \S2.2.4, the Mach configuration must be the flat one
constructed there. Theorem \ref{thm6} is then proved.

{\sc Step 2.} Now we show there will be contradictions if either
$p^{(2)}$ or $p^{(3)}$ is not constant. There are two cases:
$\Rmnum{1}.$ One is constant while the other is not; $\Rmnum{2}.$
Both are not constant.

{\sc Step 3. Case $\Rmnum{1}$.}  Without loss of generality, we
assume $p^{(2)}$ is a constant, but $p^{(3)}$ is not. Then
$p^{(2)}\equiv p_m$ in $\overline{\Omega_2}$, hence $p^{(3)}=p_m$ on
$\eta=0$.

Let $p^{(3)}_{sup}$ and $p^{(3)}_{inf}$ be respectively  the
supremum and infimum of $p^{(3)}$ in $\overline{\Omega_3}.$ If
$p^{(3)}_{sup}>p_m$ or $p^{(3)}_{inf}<p_m$, by the strong maximum
principle, they can only be attained on the shock-front
$\mathcal{S}_3\backslash\{O\}$. As shown in \S 4.3, there will be
contradictions.

{\sc Step 4. Case $\Rmnum{2}$.} Denote $p^{(2)}_{sup}$ and
$p^{(2)}_{inf}$ be respectively the supremum and infimum of
$p^{(2)}$ in $\overline{\Omega_2}.$

{\sc Step 4.1.} We firstly show $p^{(2)}_{sup}>p_m$ implies
contradictions. Suppose now $p^{(2)}_{sup}>p_m.$ Then by the strong
maximum principle, the extreme $p^{(2)}_{sup}$ can only be attained
on the shock-front $\mathcal{S}_2\backslash\{O\}$ or on the contact
line $\eta=0.$ As shown in \S 4.3, there will be contradictions if
this happens on $\mathcal{S}_2$. So it is only possible to take
place on $\eta=0.$ By \eqref{501}, we infer $p^{(3)}_{sup}>p_m.$ As
before, $p^{(3)}_{sup}$ also can only be attained on $\eta=0.$ Hence
we get
\begin{eqnarray}\label{508}
p^{(3)}_{sup}=p^{(2)}_{sup}
\end{eqnarray}
and they are attained at the same point $Q$ ($Q\ne O$) on the
contact line. Therefore, by Lemma \ref{Hopf}, we have
\begin{eqnarray}
&\p_\eta p^{(3)}(Q)<0,\label{509}\\
&\p_\eta p^{(2)}(Q)>0.\label{510}
\end{eqnarray}
Now let us look at the equations \eqref{503}\eqref{504}. We may
solve (cf. \eqref{a21})
\begin{eqnarray}
&&\p_\eta p^{(i)}=-\frac{1}{\beta_2^{(i)}}(\p_\xi
w^{(i)}+\lambda_R^{(i)}\p_\eta w^{(i)}),\\
&&\p_\xi p^{(i)}=\frac{1}{\beta_2^{(i)}}(\lambda_R^{(i)}\p_\xi
w^{(i)}+(\beta_1^{(i)}\beta_2^{(i)}+(\lambda_R^{(i)})^2)\p_\eta
w^{(i)}).
\end{eqnarray}
Then by \eqref{503}\eqref{504} and \eqref{509}\eqref{510}, we may
get, at $Q$,
\begin{eqnarray}
&(\p_\xi w^{(3)}, \p_\eta w^{(3)})\bullet(1,\lambda_R^{(3)})>0,\\
&(\p_\xi w^{(2)}, \p_\eta w^{(2)})\bullet(1,\lambda_R^{(2)})<0.
\end{eqnarray}
On the other hand, note that we have $\p_\xi p^{(2)}(Q)=\p_\xi
p^{(3)}(Q)=0$ by \eqref{508}, therefore at $Q$,
\begin{eqnarray}
(\p_\xi w^{(i)}, \p_\eta w^{(i)})\bullet(\lambda_R^{(i)},
\beta_1^{(i)}\beta_2^{(i)}+(\lambda_R^{(i)})^2)=0.
\end{eqnarray}
We then get
\begin{eqnarray}
\p_\xi w^{(i)}(Q)&=&(\p_\xi w^{(i)}, \p_\eta w^{(i)})\bullet(1,0)\nonumber\\
&=&-\frac{\lambda_R^{(i)}}{\beta_1^{(i)}\beta_2^{(i)}}(\p_\xi
w^{(i)}, \p_\eta w^{(i)})\bullet(\lambda_R^{(i)},
\beta_1^{(i)}\beta_2^{(i)}+(\lambda_R^{(i)})^2)\nonumber\\
&&+\frac{\beta_1^{(i)}\beta_2^{(i)}+(\lambda_R^{(i)})^2}{\beta_1^{(i)}\beta_2^{(i)}}
(\p_\xi w^{(i)}, \p_\eta
w^{(i)})\bullet(1,\lambda_R^{(i)})\nonumber\\
&=&\frac{\beta_1^{(i)}\beta_2^{(i)}+(\lambda_R^{(i)})^2}{\beta_1^{(i)}\beta_2^{(i)}}
(\p_\xi w^{(i)}, \p_\eta
w^{(i)})\bullet(1,\lambda_R^{(i)})\begin{cases}
<0 & i=2,\\
>0 &  i=3.
\end{cases}\nonumber\\
&& \label{516}
\end{eqnarray}
However, by \eqref{501}, we should have $\p_\xi w^{(2)}(Q)=\p_\xi
w^{(3)}(Q),$ a contradiction!

{\sc Step 4.2.} We have now $p^{(2)}_{sup}=p_m.$ Since $p^{(2)}$ is
not constant, one has $p^{(2)}_{inf}<p_m.$ Arguments as in {\sc Step
4.1} show that we must have $p^{(2)}_{inf}=p^{(3)}_{inf}$ and they
are attained at the same point $Q\ne O$ on the contact line, and a
contradiction similar to \eqref{516} then follows.

This finishes the proof of Theorem \ref{thm6}.

\bigskip
{\bf Acknowledgments.} B. Fang is supported in part by NNSF of China
under Grant No. 10801096.  L. Liu is supported in part by NNSF of
China under Grant No. 10971134.  H. Yuan is supported in part by
NNSF of China under Grants No. 10901052, No. 10871071, and Chenguang
Program (09CG20) sponsored by Shanghai Municipal Education
Commission and Shanghai Educational Development Foundation.



\begin{thebibliography}{99}
\bibitem{Ben}
Ben-Dor, G. \emph{Shock Wave Reflection Phenomena.} Second edition.
Springer, Berlin, 2007.


\bibitem{chenchenf}
 Chen, G.-Q.; Chen, J.; Feldman, M. Transonic shocks and free
 boundary problems for the full Euler equations in infinite nozzles.
\emph{ J. Math. Pures Appl.} \textbf{88} (2007), no. 2, 191--218.

\bibitem{cf2003}
 Chen, G.-Q.; Feldman, M. Multidimensional transonic
 shocks and free boundary problems for nonlinear equations of mixed type.
\emph{ J. Amer. Math. Soc.} \textbf{16} (2003), no. 3, 461--494
(electronic).


\bibitem{cf2004}
 Chen, G.-Q.; Feldman, M. Steady transonic shocks and
  free boundary problems for the Euler equations in infinite cylinders.
 \emph{Comm. Pure Appl. Math.} \textbf{57} (2004), no. 3, 310--356.


\bibitem{cf2007}
 Chen, G.-Q.; Feldman, M. Existence and stability of multidimensional
  transonic flows through an infinite nozzle of arbitrary
 cross-sections. \emph{Arch. Ration. Mech. Anal.} \textbf{184} (2007), no. 2, 185--242.


\bibitem{chenli}
Chen, G.-Q.; Li, T.-H. Well-posedness for two-dimensional steady
supersonic Euler flows past a Lipschitz wedge. \emph{J. Differential
Equations} \textbf{244} (2008), no. 6, 1521--1550.


\bibitem{CYu}
Chen, G.-Q.; Yuan, H.  Uniqueness of transonic shock solutions in a
duct for steady potential flow. {\em J. Differential Equations}
\textbf{247} (2009), no. 2, 564--573.


\bibitem{czz}
Chen, G.-Q.; Zhang, Y.; Zhu, D. Existence and stability of
supersonic Euler flows past Lipschitz wedges. \emph{Arch. Ration.
Mech. Anal.} \textbf{181} (2006), no. 2, 261--310.

\bibitem{ChenJun}
Chen, J. Subsonic flows for the full Euler equations in half plane.
\emph{J. Hyperbolic Differ. Equ.} \textbf{6} (2009), no. 2,
207--228.

\bibitem{chens1}
Chen, S. Stability of transonic shock fronts in two-dimensional
Euler systems. \emph{Trans. Amer. Math. Soc.} \textbf{357} (2005),
no. 1, 287--308 (electronic).

\bibitem{chens2}
Chen, S. Stability of a Mach configuration. \emph{Comm. Pure Appl.
Math.} \textbf{59} (2006), no. 1, 1--35.


\bibitem{chen3}
Chen, S. Mach configuration in pseudo-stationary compressible flow.
\emph{J. Amer. Math. Soc.} \textbf{21} (2008), no. 1, 63--100
(electronic).



\bibitem{CHEN-FANG}
Chen, S.; Fang, B. Stability of transonic shocks in supersonic flow
past a wedge. \emph{J. Differential Equations} \textbf{233} (2007),
no. 1, 105--135.



\bibitem{ChenFang}
Chen, S.; Fang, B. Stability of reflection and refraction of shocks
on interface. {\it J. Differential Equations} \textbf{244} (2008),
no. 8, 1946--1984.


\bibitem{CF}
Courant, R.; Friedrichs, K.~O. {\it Supersonic Flow and Shock
Waves}. Interscience Publishers,  New York, 1948.


\bibitem{Da}
 Dafermos, C.~M. {\em Hyperbolic Conservation Laws in Continuum
 Physics}. Springer-Verlag, New York, 2000.


\bibitem{Fang}
Fang, B. {Stability of transonic shocks for the full Euler equations
in supersonic flow past a wedge}. \emph{Math. Methods Appl. Sci.}
\textbf{29} (2006), no. 1, 1--26.


\bibitem{FG}
Finn, R.; Gilbarg, D. Asymptotic behavior and uniqueness of plane
subsonic flows. \emph{Comm. Pure Appl. Math.} \textbf{10}(1957),
23--63.

\bibitem{GT}
Gilbarg, D.; Trudinger, N.~S. {\em Elliptic Partial Differential
Equations of Second Order}. Reprint of the 1998 edition. Springer,
Berlin, 2001.


\bibitem{Kim}
Kim, E. H. Existence and stability of perturbed transonic shocks for
compressible steady potential flows. \emph{Nonlinear Anal.}
\textbf{69} (2008), no. 5-6, 1686--1698.


\bibitem{LIYU}
{Li, T.-T.;  Yu, W.-C. } {\it Boundary value problems for
quasilinear hyperbolic systems}, Duke University Mathematics Series
5. Duke University, Mathematics Department, Durham, N.C., 1985.


\bibitem{L3}
 Liu, L. Unique subsonic compressible potential flows in
three--dimensional ducts, \emph{Discrete and Continuous Dynamical
Systems} \textbf{27} (2010), no. 1, 357--368.


\bibitem{LY2}
{ Liu, L.;  Yuan, H.}  {Uniqueness of  Symmetric Steady Subsonic
Flows in Infinitely Long Divergent Nozzles},  (2009), submitted.


\bibitem{Wjh}
Wang, J.-H. {\it Two--dimensional unsteady flows and shock
 waves}. Scince Press, Beijing, 1994 (in Chinese).


\bibitem{XYY2009}
 Xin, Z.; Yan, W.; Yin, H. Transonic shock
 problem for the Euler system in a nozzle.
\emph{ Arch. Ration. Mech. Anal.} \textbf{194} (2009), no. 1, 1--47.


\bibitem{XY2005}
 Xin, Z.; Yin, H. Transonic shock in a nozzle.
 I. Two-dimensional case. \emph{Comm. Pure Appl. Math.} \textbf{58} (2005),
 no. 8, 999--1050.


\bibitem{XY2008}
 Xin, Z.; Yin, H.  The transonic shock in a nozzle,
 2-D and 3-D complete Euler systems.
\emph{ J. Differential Equations} \textbf{245} (2008), no. 4,
1014--1085.


\bibitem{Yin}
Yin, H.; Zhou, C.  On global transonic shocks for the steady
supersonic Euler flows past sharp 2-D wedges. \emph{J. Differential
Equations }\textbf{246} (2009), no. 11, 4466--4496.


\bibitem{Yu2}
 Yuan, H. {On transonic shocks in two--dimensional variable--area
ducts for steady Euler system}, {\it SIAM J. Math. Anal.} \textbf{
38} (2006), 1343--1370.


\bibitem{Zhangyongqian2003}
Zhang, Y. Steady supersonic flow past an almost straight wedge with
large vertex angle. \emph{J. Differential Equations} \textbf{192}
(2003), no. 1, 1--46.


\end{thebibliography}
\end{document}